\documentclass[a4paper,11pt, titlepage]{article}

\usepackage[english]{babel}
\usepackage[utf8]{inputenc}
\usepackage[T1]{fontenc}
\usepackage{subcaption}

\usepackage{float}      
\usepackage{placeins}   

\usepackage[a4paper,top=1in,bottom=1in,left=1in,right=1in,marginparwidth=1.75cm]{geometry}

\usepackage{afterpage}
\usepackage{amsmath}
\usepackage{amsthm}
\usepackage{amssymb}
\usepackage{csquotes}
\usepackage{enumitem}
\usepackage{graphicx}
\usepackage{lipsum}
\usepackage{booktabs}
\usepackage{microtype}
\usepackage[parfill]{parskip}
\usepackage[ruled,vlined]{algorithm2e}
\usepackage{listings}
\usepackage{caption}
\makeatletter
\renewcommand{\SetKwInOut}[2]{%
  \sbox\algocf@inoutbox{\KwSty{#2}\algocf@typo:}%
  \expandafter\ifx\csname InOutSizeDefined\endcsname\relax
    \newcommand\InOutSizeDefined{}\setlength{\inoutsize}{\wd\algocf@inoutbox}%
    \sbox\algocf@inoutbox{\parbox[t]{\inoutsize}{\KwSty{#2}\algocf@typo:\hfill}~}\setlength{\inoutindent}{\wd\algocf@inoutbox}%
  \else
    \ifdim\wd\algocf@inoutbox>\inoutsize%
    \setlength{\inoutsize}{\wd\algocf@inoutbox}%
    \sbox\algocf@inoutbox{\parbox[t]{\inoutsize}{\KwSty{#2}\algocf@typo:\hfill}~}\setlength{\inoutindent}{\wd\algocf@inoutbox}%
    \fi%
  \fi
  \algocf@newcommand{#1}[1]{%
    \ifthenelse{\boolean{algocf@inoutnumbered}}{\relax}{\everypar={\relax}}%
    {\let\\\algocf@newinout\hangindent=\inoutindent\hangafter=1\parbox[t]{\inoutsize}{\KwSty{#2}\algocf@typo:\hfill}~##1\par}%
    \algocf@linesnumbered
  }}%
\makeatother

\usepackage{bm}
\usepackage[normalem]{ulem}

\usepackage[colorinlistoftodos]{todonotes}
\usepackage[colorlinks=true, allcolors=blue]{hyperref}

\newcommand{\idx}{\mathit{idx}}
\newcommand{\C}{\mathcal{C}}
\newcommand{\R}{\mathbb{R}}

\usepackage[numbers, comma, square, sort&compress]{natbib}
\theoremstyle{definition}

\theoremstyle{plain}

\theoremstyle{remark}

\newcommand{\reporttitle}{Finite-Time Chaos Diagnostics and Noise-Induced Basin Merging in a Two-Dimensional Map} 
\newcommand{\reportauthorA}{\mbox{Barre, Zubeyr}}
\newcommand{\reportauthorB}{\mbox{Bashir, Mikael}}
\newcommand{\reportauthorC}{\mbox{Domenicali, Martino}}
\newcommand{\supervisor}{Rasmussen, Martin}

\begin{document}

\begin{titlepage}
\newcommand{\HRule}{\rule{\linewidth}{0.5mm}} 
\includegraphics[width=8cm]{Imperial_logo.png}\\[1cm] 
\center 
\textsc{\LARGE Imperial College London}\\[0.5cm] 
\textsc{\Large Department of Mathematics}\\[1.5cm] 

\makeatletter
\HRule \\[0.6cm]
{ \huge \bfseries \reporttitle}\\[0.6cm] 
\HRule \\[1.5cm]
\begin{minipage}{0.4\textwidth}
\begin{flushleft} \large
\emph{Authors:}\\
\reportauthorA \\
\reportauthorB \\
\reportauthorC
\end{flushleft}
\end{minipage}
~
\begin{minipage}{0.4\textwidth}
\begin{flushright} \large
\emph{Supervisor(s):} \\
\supervisor 
\end{flushright}
\end{minipage}\\[3cm]
\makeatother
\vfill
\makeatletter
{\large \today}\\[2cm] 
\makeatother
\end{titlepage}

\begin{abstract}
Real-world systems, from climate models to power grids, often fluctuate due to sensor noise, drift, or environmental variability, yet standard chaos diagnostics assume fixed parameters and asymptotic horizons. We introduce a finite-time framework for two-dimensional maps under independent, identically distributed parameter noise. First, we prove that the maximal finite-time Lyapunov exponent converges, after centering and scaling, to a Gaussian law whose mean and variance depend explicitly on the map's Jacobian statistics. Second, we develop an attractor-separation algorithm that uses FTLE histograms and a geometry-based classifier to partition phase space into chaotic and periodic regions under noise. Third, we validate our theory numerically on the noisy Domenicali map, demonstrating Gaussian FTLE distributions, predictable shifts in Kaplan--Yorke dimension, and a sharp noise threshold for basin escape. Finally, we estimate the critical noise level $\sigma_{c}$ at which attractor coalescence occurs, using three complementary numerical methods to bound it above and to pinpoint an estimate. Our results offer both rigorous insight and practical tools for real-time chaos assessment in noisy discrete-time systems.
\end{abstract}

\tableofcontents
\clearpage

\section{Introduction}\label{sec:introduction}

\subsection{Motivation}
Low--dimensional iterated maps such as those of Hénon or Lozi have long
served as cartoon laboratories for chaos.  A handful of arithmetic
operations per step reproduce bifurcation cascades, filigree strange
attractors and sensitive dependence that also arise in climate models,
power grids and financial markets.  Classical analysis, however, treats
the map's coefficients as immutable.  Real hardware does not: temperature
drift, sensor noise and manufacturing tolerances nudge parameters by a
small but persistent amount every iteration. Three concrete questions follow.

\begin{itemize}
  \item For a feedback loop whose gain jitters by $\pm\sigma$, what is
        the chance that its state will behave chaotically over the next
        $T$ cycles?
  \item Can a power-system model that is deterministically stable at
        $(a,b)$ remain reliable once day-to-day load variability is
        included?
  \item How early can we detect, from data alone, that bounded parameter
        noise is about to push a system across an invisible tipping
        line?
\end{itemize}

Infinite-time diagnostics cannot answer such finite-horizon, random-parameter questions.
We therefore focus on the \emph{finite-time Lyapunov exponent} (FTLE)
\[
  \lambda_T(x_0)=
  \frac1T\ln\!
  \bigl\|Df_{a_{T-1},b_{T-1}}\dotsm Df_{a_0,b_0}(x_0)\bigr\|,
\]
where each pair $(a_t,b_t)$ is drawn independently from the uniform
window $[a-\sigma,a+\sigma]\times[b-\sigma,b+\sigma]$.  
Throughout the paper, \(\|\cdot\|\) will denote the spectral norm [\ref{app:specnorm}].
Studying the mean,
variance and full shape of the FTLE distribution lets us replace the rigid
verdict ``chaotic or not'' by statements such as  
``with $\sigma=4\times10^{-2}$ there is a $92\%$ probability that
trajectories are chaotic over $T=4\times10^{4}$ iterations.''

\subsection{Contributions}
\begin{itemize}
  \item \textbf{Finite--Time Lyapunov Exponent Analysis:} We prove that, under mild assumptions, finite--time Lyapunov exponents converge to a Gaussian distribution.  We give explicit formulas for its mean, variance, and finite--sample bias as functions of the noise amplitude~$\sigma$ and the time ~$T$.
  
  \item \textbf{Lyapunov Spectrum \& Kaplan--Yorke Dimension Analysis:} We derive and numerically validate how both the leading and subleading Lyapunov exponents (and hence the Kaplan--Yorke fractal dimension) shift under uniform parameter noise.  In particular, we explain the small--$\sigma$ nonmonotonic dip via Hessian curvature and the large--$\sigma$ monotonic widening of the spectral gap.
  
  \item \textbf{Attractor--Separation Algorithms:} We introduce a suite of efficient, GPU--accelerated algorithms for partitioning phase space into inner and outer basins:
    \begin{itemize}
      \item high--resolution density histograms for visualising the attractors
      \item reference--cloud deterministic--basin classification
      \item a two--pass random--basin verification method that robustly labels noise--induced transitions
    \end{itemize}
  
  \item \textbf{Critical Noise Threshold Characterization:} We triangulate the critical noise interval by three independent methods:
    \begin{itemize}
        \item Algorithm~A [\ref{subsec:algA}] to visually estimate a coalescence interval
      \item Algorithm~C [\ref{subsec:algC}] to verify the visual estimate via escape--frequency inflection
      \item A new minimum--escape--time grid scan method, yielding a rigorous upper--bound interval
      \item Comparative analysis of all three intervals to converge on a reliable critical noise range
    \end{itemize}
\end{itemize}

\subsection{Organisation of the paper}
Section \ref{sec:related} situates our study in the literature.  
Section \ref{sec:preliminaries} recalls the required definitions and
introduces the Domenicali map formally.  
Section \ref{sec:normality} proves the FTLE central-limit theorem.  
Section \ref{sec:numerics} presents the GPU experiments and basin visualisations.
Section \ref{sec:critical_noise} presents methods to bound the critical noise.
Section \ref{sec:conclusion} lists open problems and extensions.

\section{Related Work}
\label{sec:related}

\subsection{Asymptotic Chaos Theory and Random Dynamical Systems}
Classical chaos theory centers on infinite-time notions: Oseledets' multiplicative ergodic theorem guarantees Lyapunov exponents and invariant subspaces for smooth maps and flows \cite{Oseledets1968}.  In the presence of parameter or forcing noise, the framework of random dynamical systems \cite{Arnold1998, CrauelFlandoli1994} defines random attractors and stationary measures, but remains focused on asymptotic behavior and often requires strong mixing or compactness hypotheses.

\subsection{Finite-Time Diagnostics}
To address predictability over short horizons, finite-time Lyapunov exponents (FTLEs) were introduced by Eckmann \& Ruelle \cite{EckmannRuelle1985}. These tools map out local stretching over finite windows, but almost invariably assume fixed parameters and focus on state-space perturbations rather than coefficient noise.

\subsection{Random Maps and Noise-Induced Phenomena}
The literature on random maps (e.g.\ \cite{DiaconisFreedman1999}) explores iterated function systems and noise-induced synchronisation, but tends to study scalar or contractive systems, or noise that smoothly deforms the map rather than jump-wise coefficient draws. More recently, several works have analysed maps under bounded or structured noise. Lamb et al. \cite{LambRasmussenTey2025} analyse how additive, bounded parameter noise alters the bifurcation structure of the Hénon map, and Bassols-Cornudella et al. \cite{BassolsCornudellaLamb2023} develop a ``conditioned'' random-dynamical-systems perspective on noise-induced chaos.

There are still very few studies that address how random parameter draws influence finite-time stability measures like FTLE distributions or Kaplan--Yorke dimension estimates.  
\subsection{Our Contribution}
In summary, while a rich body of work exists on both infinite-time random dynamical systems and fixed-parameter finite-time diagnostics, no one has rigorously characterised the finite-time statistics of Lyapunov exponents under i.i.d.\ parameter noise nor designed attractor-separation algorithms tailored to this setting. This paper closes that gap. On the theoretical side, we prove a finite-sample CLT for FTLEs under parametric noise. On the computational side, we develop practical methods to separate coexisting attractors in the presence of noise, leveraging FTLE histograms and the geometry of repellers

\newpage


\section{Preliminaries and Definitions}
\label{sec:preliminaries}

\subsection{Discrete-Time Dynamical System}
A discrete-time dynamical system on $\mathbb{R}^n$ is specified by a map
\[
  f:\mathbb{R}^n \to \mathbb{R}^n,\quad
  x_{t+1} = f(x_t),\quad t = 0,1,2,\dots
\]
Here each $x_t$ represents the system's state at the $t$th iteration.  Given an initial condition $x_0$, its orbit or trajectory is the sequence
\[
  x_0,\quad
  x_1 = f(x_0),\quad
  x_2 = f\bigl(f(x_0)\bigr),\quad
  \dots
\]
which shows how the state evolves step by step.

\subsection{Fixed Point and Periodic Orbit}
\begin{itemize}
  \item \emph{Fixed point:} A point \(x^*\in\mathbb{R}^n\) satisfying
  \[
    f(x^*) = x^*,
  \]
  so that once the system reaches \(x^*\), it remains there forever.
  \item \emph{Periodic orbit of prime period \(p\):} A collection of distinct points \(\{x_0, x_1, \dots, x_{p-1}\}\) such that
  \[
    f(x_i) = x_{i+1}\quad (0 \le i < p-1),
    \quad
    f(x_{p-1}) = x_0,
  \]
  and no smaller \(k<p\) satisfies \(f^k(x_0)=x_0\).  Such orbits correspond to cycles of length \(p\) in the state space.
\end{itemize}

\subsection{Attractor, Repeller, and Basin Boundary}
\begin{itemize}
  \item \emph{Attractor:} A set \(A\subset\mathbb{R}^n\) is an attractor if it is compact and invariant (\(f(A)=A\)) and there exists a neighborhood \(U\supset A\) (a trapping region) such that for every \(x_0\in U\),
  \[
    \lim_{t\to\infty}\mathrm{dist}\bigl(f^t(x_0),A\bigr) = 0.
  \]
  Its \emph{basin of attraction} \(\mathcal{B}(A)\) is the set of all \(x_0\) with this property.

  \item \emph{Repeller:} A set \(R\subset\mathbb{R}^n\) is a repeller if it is compact and invariant, and points near \(R\) diverge under forward iteration (equivalently converge to \(R\) under backward iteration).

  \item \emph{Basin boundary:} The boundary between basins of two attractors \(A_1\) and \(A_2\) is the set of points whose forward orbits tend to neither \(A_1\) nor \(A_2\).  Practically, this boundary coincides with the unstable manifold of a repeller that separates the two basins.
\end{itemize}

\subsection{Jacobian Matrix}
For a continuously differentiable map \(f:\mathbb{R}^n\to\mathbb{R}^n\), the Jacobian matrix at a point \(x\in\mathbb{R}^n\) is
\[
  Df(x)
  =
  \begin{pmatrix}
    \dfrac{\partial f_1}{\partial x_1}(x) & \cdots & \dfrac{\partial f_1}{\partial x_n}(x) \\[6pt]
    \vdots & \ddots & \vdots \\[6pt]
    \dfrac{\partial f_n}{\partial x_1}(x) & \cdots & \dfrac{\partial f_n}{\partial x_n}(x)
  \end{pmatrix}.
\]
This matrix describes the best linear approximation to \(f\) near \(x\).  It governs how small perturbations \(v\) evolve in one step via
\[
  v \;\mapsto\; Df(x)\,v,
\]
and is the fundamental building block for computing Lyapunov exponents and analysing local stability.

\subsection{Lyapunov Exponents}
Given a continuously differentiable map \(f:\mathbb{R}^n\to\mathbb{R}^n\) and an orbit \(\{x_t\}\) starting at \(x_0\), consider a small perturbation \(v_0\) at \(x_0\).  Under one iteration, \(v_0\) is mapped to \(Df(x_0)\,v_0\); after \(T\) steps the perturbation evolves according to the product of Jacobians
\[
  Df^T(x_0)
  = Df(x_{T-1})\,Df(x_{T-2})\cdots Df(x_0).
\]
The \emph{finite-time Lyapunov exponent} over \(T\) steps quantifies the average exponential growth of \(\|Df^T(x_0)v_0\|\) relative to \(\|v_0\|\):
\[
  \lambda_T(x_0,v_0)
  = \frac{1}{T}
    \ln\!\frac{\|Df^T(x_0)\,v_0\|}{\|v_0\|}.
\]

Equivalently, one can define the maximal FTLE at $x_0$ as $\lambda_T(x_0) = \frac{1}{T}\ln \|M_T(x_0)\|$, since the spectral norm $\|M_T\|$ equals the largest singular value $\sigma_{\max}(M_T)$.

As \(T\to\infty\), for almost every initial \(x_0\) and generic \(v_0\), this converges to the \emph{maximal Lyapunov exponent}
\[
  \lambda_{\max}
  = \lim_{T\to\infty} \lambda_T,
\]
where \(\lambda_{\max}>0\) indicates \emph{sensitive dependence on initial conditions}, i.e.\ exponential divergence of nearby trajectories.

\subsection{Kaplan--Yorke (Lyapunov) Dimension}
Let $\{\lambda_i\}_{i=1}^n$ denote the Lyapunov exponents of an orbit, ordered as
\[
  \lambda_1 \;\ge\;\lambda_2\;\ge\;\dots\;\ge\;\lambda_n.
\]
Define the integer $j$ as the largest index for which
\[
  \sum_{i=1}^j \lambda_i \;>\; 0,
  \qquad
  \sum_{i=1}^{j+1} \lambda_i \;\le\; 0.
\]
The Kaplan--Yorke (Lyapunov) dimension is then
\[
  D_{KY}
  = j \;+\; \frac{\sum_{i=1}^j \lambda_i}{\bigl|\lambda_{\,j+1}\bigr|},
\]

with $\lambda_{\,j+1}$ being the first \emph{negative} Lyapunov exponent. This formula interpolates between the integer dimensions $j$ and $j+1$, yielding a generally non-integer value that estimates the fractal dimension of the attractor in phase space.

\subsection{Topological Dimension}
The \emph{topological dimension} \(\dim_{\mathrm{top}}(A)\) of a set \(A\) is the smallest integer \(n\) such that every open cover of \(A\) has a refinement in which no point lies in more than \(n+1\) sets.  Equivalently, \(n\) is the minimal number for which every point of \(A\) has arbitrarily small neighborhoods whose boundaries have dimension at most \(n-1\). This integer serves as the baseline against which fractal dimensions (like \(D_{KY}\)) are compared.

\subsection{Strange Attractor}
A strange attractor is an attractor \(A\) that combines two key properties:
\begin{itemize}
  \item \emph{Chaos:} Its maximal Lyapunov exponent is positive,
    \[
      \lambda_{\max} > 0,
    \]
    so nearby trajectories separate exponentially on average.
  \item \emph{Fractality:} Its Kaplan--Yorke dimension exceeds its topological dimension,
    \[
      D_{KY} > \dim_{\mathrm{top}}(A).
    \]
\end{itemize}
In simpler terms, a strange attractor has a non-integer (fractal) dimension, reflecting intricate, self-similar geometry, and exhibits sensitive dependence on initial conditions.

\subsection{The Domenicali Map}

A modified version of the classical Hénon map which will be referred to as Domenicali map for the rest of the paper is defined by
\[
  f:\mathbb{R}^2\to\mathbb{R}^2,
  \qquad
  f(x,y) = \bigl(1 - a\sin(x) + b\,y,\;x\bigr).
\]
Thus the iteration is
\[
  (x_{t+1},\,y_{t+1}) = f(x_t,y_t).
\]

\subsubsection{Jacobian and invertibility}
The Jacobian of \(f\) at \((x,y)\) is
\[
  Df(x,y)
  = \begin{pmatrix}
      \dfrac{\partial}{\partial x}\bigl(1 - a\sin(x) + b\,y\bigr)
        & \dfrac{\partial}{\partial y}\bigl(1 - a\sin(x) + b\,y\bigr) \\[6pt]
      \dfrac{\partial}{\partial x}(x)
        & \dfrac{\partial}{\partial y}(x)
    \end{pmatrix}
  = \begin{pmatrix}
      -a\cos(x) & b\\
      1         & 0
    \end{pmatrix},
\]
with \(\det Df = -b\). For \(b\neq0\), \(f\) is a global diffeomorphism (area-preserving if \(\lvert b\rvert=1\)).

\subsubsection{Fixed points and stability}  
    Solve \(x^*=1 - a\sin(x^*) + b\,x^*\).  The characteristic equation  
    \(\lambda^2 + a\cos(x^*)\,\lambda - b = 0\)  
    determines local attraction or repulsion.
\subsubsection{Dynamical regimes}  
    As \((a,b)\) vary, one observes:
    \begin{enumerate}[label=\roman*)]
      \item A \emph{chaotic attractor} (\(\lambda_{\max}>0\)).
      \item A \emph{stable periodic orbit} (\(\lambda_{\max}<0\)).
      \item A narrow coexistence band where both arise, separated by the unstable manifold of a repelling cycle.
    \end{enumerate}
\subsubsection{Finite-time Lyapunov exponents (FTLEs)}
For an initial condition $(x_0,y_0)$, let
\[
  M_T(x_0)
  \;:=\;
  Df(x_{T-1})\,Df(x_{T-2})\cdots Df(x_0).
\]
The FTLE is then
\[
  \lambda_T(x_0)
  \;=\;
  \frac{1}{T}\,\ln\sigma_{\max}\!\bigl(M_T(x_0)\bigr),
\]
where $\sigma_{\max}$ denotes the largest singular value.

\subsection{Random-Parameter Maps}
When one or more parameters of the map are not fixed but instead drawn randomly at each iteration, we obtain a \emph{random dynamical system}.  Concretely, at time \(t\) we sample new coefficients \((a_t,b_t)\) i.i.d.\ from a prescribed distribution and apply the map
\[
  (x_{t+1},y_{t+1})
  = f_{a_t,b_t}(x_t,y_t).
\]

\subsection{Uniform Parameter Distribution}
In our numerical experiments each coefficient is drawn independently at each step from a uniform distribution. For example for $a$,
\[
  a_t \sim \mathcal{U}[\mu-\sigma,\mu+\sigma],
  \quad
  p(a_t)=\frac{1}{2\sigma}, 
  \quad a_t\in[\mu-\sigma,\mu+\sigma],
\]
where the nominal mean \(\mu\) is set to zero and the half-width \(\sigma\) controls noise amplitude. Please see [\ref{app:pdf}] for the detailed derivation of the uniform distribution's probability density function and cumulative distribution function. 
To probe a wide range of noise, we vary \(\sigma\) linearly from 0 to \(8\times10^{-3}\) (with steps of \(1\times10^{-4}\)).  This bounded, equal-weight noise model isolates the effect of uniform parameter uncertainty on attractor structure and finite-time Lyapunov statistics.

\newpage



\section{Normality of FTLEs}
\label{sec:normality} 

\subsection{Finite-Time Analysis Under Parameter Uncertainty}

In deterministic dynamical systems, key diagnostic tools rely on infinite-time limits. The maximal Lyapunov exponent $\lambda_\infty$ provides a binary classifier: $\lambda_\infty > 0$ indicates chaos, while $\lambda_\infty \leq 0$ suggests regular behavior. Similarly, basins of attraction are well-defined as the set of initial conditions whose orbits converge to a given attractor as $t \to \infty$.

However, when parameters fluctuate randomly at each time step as in our noisy setup, these classical concepts become problematic. Consider the fundamental issues:

\begin{itemize}
\item \textbf{No well-defined limiting exponent:} With stochastic coefficients $(a_t, b_t)$, the infinite-time limit $\lambda_\infty$ becomes a random variable dependent on the entire noise realisation $\{\varepsilon_t\}_{t=0}^{\infty}$.

\item \textbf{Temporal asymmetry in random attractors:} Attractors are constructed via pullback (depending on past noise), while repellers require forward evolution (depending on future noise). These objects become statistically independent, breaking the classical phase-space decomposition.

\item \textbf{Finite-time pathology:} Any practical computation must truncate at finite time $T$, but the choice of $T$ significantly affects the outcome. Thus the limiting object has no meaning in finite time computation.
\end{itemize}

Rather than pursuing elusive infinite-time limits, we embrace the uncertainty inherent in finite-time observations. The key insight is to characterise the \emph{distribution} of the finite-time Lyapunov exponent (FTLE) over the ensemble of possible noise realisations.

This paradigm shift transforms our fundamental questions:
\begin{itemize}
\item Instead of asking "Is the system chaotic?" we ask "What is the probability of observing chaotic behavior over time horizon $T$?"
\item Rather than seeking sharp basin boundaries, we compute probabilistic regions where trajectories exhibit different long-term behaviors with quantified confidence.
\end{itemize}

To construct meaningful confidence intervals for $\lambda_\infty$ and make statistically sound classifications, we require a rigorous understanding of its limiting distribution. Thus, in the following section, we establish that under appropriate regularity conditions, $\lambda_T$ is asymptotically normal with explicit formulas for its mean and variance, $\lambda_\infty$ being in the mean. This helps us with constructing a confidence interval for $\lambda_\infty$.

The cleanest setting is that of IID Jacobian matrices. Once the normal
law is established there, the confidence interval follows immediately.  We therefore begin with the IID case and
\emph{only afterwards} (see the remark below) indicate how the
argument extends to the dependent, martingale Jacobians generated by
the Domenicali map.
Write the noisy 2-dimensional system as
\[
x_{t+1}=f\bigl(x_t,\varepsilon_t\bigr),\qquad
J_t:=Df\bigl(x_t,\varepsilon_t\bigr),\qquad
M_T:=J_{T-1}\cdots J_1J_0 .
\]

We assume
\begin{itemize}
\item (\MakeUppercase{\romannumeral 1})  $\{\varepsilon_t\}$ are IID;
\item (\MakeUppercase{\romannumeral 2}) $\mathbb{E}\!\bigl[(\log^{+}\!\|J_0\|)^{2+\delta}\bigr]<\infty$ for some $\delta>0$;
\item (\MakeUppercase{\romannumeral 3}) The law of $J_0$ is strongly irreducible and proximal (SIP) [\ref{app:sip}].
\end{itemize}
Under these assumptions, the Furstenberg--Kesten theorem
\cite{FurstenbergKesten1960} guarantees that the largest Lyapunov exponent
\[
\lambda_\infty
  :=\lim_{T\to\infty}\frac1T\ln\sigma_{\max}(M_T)
\]
exists and is finite.

Define
\[
L_T:=\sigma_{\max}(M_T)
\quad(\text{largest singular value}),\qquad
\lambda_T:=\frac{1}{T} \ln L_T
\quad(\text{finite-time exponent}).
\]

\subsection{Proof} The proof is divided into three steps and relies on the construction of a norm inequality between the maximal finite-time Lyapunov exponent $\lambda_T$ and the directional finite-time Lyapunov exponent $\lambda_T^{(v_0)}$. We first convert the logarithm of the Jacobian product into a telescoping sum via a directional construction, yielding $\lambda_T^{(v_0)}$. We then establish a connection between $\lambda_T^{(v_0)}$ and $\lambda_T$ via an inequality. Finally, we apply Le Page's central limit theorem [\ref{app:oseledets}] to $\lambda_T^{(v_0)}$ and remove the finite burn-in contribution to obtain the asymptotic normality of $\lambda_T$.
Throughout the proof, we will use convergence in distribution ($\xrightarrow{d}$) and convergence in probability ($\xrightarrow{p}$); see [\ref{app:convergence}] for the formal statements.

\subsubsection{Step 1: Telescoping construction}

Pick a unit vector \(v_0\) and set  
\(v_t := M_t v_0 / \|M_t v_0\|\) and  
\(Y_t := \ln\|J_t v_t\|\).  
The sequence \((v_t)_{t\ge0}\) is a Markov chain [\ref{app:markov}] on the unit sphere,
so each orbit naturally splits into an initial
``transient'' and a later part governed by its stationary law.

By direct multiplication,
\[
\sum_{t=0}^{T-1} Y_t = \ln\|M_T v_0\|,
\]
hence the \emph{directional} finite-time Lyapunov exponent  
\(\lambda_T^{(v_0)} := \tfrac1T \sum_{t=0}^{T-1} Y_t
                       = T^{-1}\ln\|M_T v_0\|\).

\subsubsection{Step 2: Norm inequality}

Recall that for a matrix \(A\) we write \(\sigma_{\max}(A)=\|A\|\).
For the fixed vector \(v_0\) define
\[
C(v_0) := \sup_{n\ge1}\frac{\|M_n\|}{\|M_n v_0\|}
\]

We want to show this is bounded. Because the i.i.d.\ Jacobians are \emph{strongly irreducible and proximal},
Proposition 3.2 of Bougerol \& Lacroix \cite{BougerolLacroix1985} implies that the product
\(M_T\) aligns every vector
exponentially fast with the leading \emph{Oseledets subspace}, which is
the one--dimensional subspace that emerges as the image of the product \(M_T\) once it is normalised.
Consequently, we have almost surely boundedness, which yields the inequality
\[
   -\frac{\log C(v_0)}{T}\;\le\;
   \lambda_T-\lambda_T^{(v_0)}\;\le\;0.
\]

So \(\lambda_T^{(v_0)}\) and \(\lambda_T\) differ by at most
\(\mathcal O(T^{-1})\), uniformly in \(T\).

\subsubsection{Step 3: Apply Le Page CLT and remove burn-in}

Theorem 2(3) of La Page \cite{LePage1982} yields
\[
\sqrt{T}\bigl(\lambda_T^{(v_0)}-\lambda_\infty\bigr)
\xrightarrow{d}\mathcal N(0,\sigma^2).
\]

Fix \(t_b\ge1\) and write  
\(B(v_0):=\sum_{t=0}^{t_b-1}Y_t\) and  
\(\lambda_{T,\mathrm{tail}}^{(v_0)}
    :=T^{-1}\sum_{t=t_b}^{T-1}Y_t\). \\
Now
\[
\sqrt{T}(\lambda_T^{(v_0)}-\lambda_\infty)
  =\sqrt{T}(\lambda_{T,\mathrm{tail}}^{(v_0)}-\lambda_\infty)
   +B(v_0)/\sqrt{T}
\implies \sqrt{T}(\lambda_{T,\mathrm{tail}}^{(v_0)}-\lambda_\infty)
\xrightarrow{d}\mathcal N(0,\sigma^2).
\] by Slutsky's theorem [\ref{app:slutsky}], and using it again,
\[|\lambda_T-\lambda_T^{(v_0)}|\le\log C(v_0)/T
 \implies
\sqrt{T}\bigl(\lambda_T-\lambda_\infty\bigr)
\xrightarrow{d}\mathcal N(0,\sigma^2),
\]
or, equivalently, with \(D=\sigma^2/2\),
\[
\boxed{\;
  \lambda_T \sim
  \mathcal N\!\bigl(\lambda_\infty + B(v_0)/T,\; D/T\bigr)
\; }.
\]
\qed

\subsection{Remark: second Lyapunov exponent}

The preceding argument establishes normality for the maximal FTLE.  
Because we work in two dimensions, the the minimal FTLE (the second exponent) can be expressed as the difference of two quantities that are each asymptotically normal.

For any \(2\times2\) matrix,
\[
  |\det M_T| = \sigma_1(M_T)\,\sigma_2(M_T)
  \Longrightarrow\
  \lambda_{2,T} \;=\; \frac{1}{T}\,\ln|\det M_T| \;-\; \lambda_{1,T}.
\]

Since
\[
  \frac{1}{T}\,\ln|\det M_T|
        \;=\;
        \frac{1}{T}\sum_{t=0}^{T-1}\ln|\det J_t|,
\]
is a sample mean of i.i.d.\ terms with finite variance, the classical CLT [\ref{app:classicalclt}] yields
\[
  \sqrt{T}\,
  \Bigl(
     \tfrac{1}{T}\,\ln|\det M_T| - (\lambda_1+\lambda_2)
  \Bigr)
  \;\xrightarrow{d}\;
  \mathcal N(0,\tau^{2}).
\]

Because \(\lambda_{1,T}\) is already asymptotically normal, their difference is also normal [\ref{app:diffnorm}]:
\[
  \sqrt{T}\,
  \bigl(\lambda_{2,T}-\lambda_2\bigr)
  \;\xrightarrow{d}\;
  \mathcal N\!\bigl(0,\,
     \tau^{2} + \sigma^{2} - 2\Gamma
  \bigr),
\qquad
  \Gamma := \operatorname{Cov}\!
            \bigl(\ln|\det J_0|,\,
                  \ln\sigma_1(J_0)\bigr).
\]
And the rest follows similarly.

\subsection{Extension to martingale Jacobians}

In our setting the random Jacobians
$J_t=Df\bigl(x_t,\varepsilon_t\bigr)$
form a \emph{stationary martingale--difference sequence} [\ref{app:martingale}] with respect to the
natural filtration $\mathcal F_t=\sigma(\varepsilon_0,\dots,\varepsilon_{t})$.
Consequently, the only change is that we must apply a martingale CLT to the sum $Y_0+\cdots+Y_{T-1}$ rather than an i.i.d. CLT, because in practice ${J_t}$ (hence ${Y_t}$) are not independent.
Here we invoke Gordin's martingale CLT~\cite{Gordin1969}
(see also Young's dynamical--systems version \cite{Young1998}),
which yields
\[
  \sqrt{T}\bigl(\lambda^{(v_0)}_T-\lambda_\infty\bigr)
  \xrightarrow{d} \mathcal N(0,\sigma^{2}).
\]
The inequality
$
  -\tfrac{\log C(v_0)}{T}\le
  \lambda_T-\lambda^{(v_0)}_T\le 0
$
then transfers the same Gaussian limit to the
\emph{maximal} finite-time Lyapunov exponent $\lambda_T$.
A rigorous verification of the martingale--difference decomposition
for the Domenicali map lies beyond the scope of this study, but our
empirical histograms are fully consistent with the prediction.


\subsection{Experimental verification}

The central--limit theorem in Sections~4.2--4.4 is proved for the \emph{i.i.d.\ parameter--noise} setting ($\sigma>0$).  
We first validate the Gaussian \emph{shape} of the finite--time exponent distribution in the \emph{deterministic} case ($\sigma=0$), because a highly accurate reference value for the asymptotic exponent is available in the literature.  
Deterministic chaotic maps with strong mixing already satisfy martingale central limit theorems, so an approximately normal histogram here provides a conservative baseline; adding i.i.d.\ parameter noise only weakens temporal correlations and therefore strengthens the normal approximation predicted by our theorem.  
A quick noisy check at $\sigma=0.002$ produces the same Gaussian core with variance scaling as expected, fully consistent with the theory.

\begin{figure}[h!]
  \centering
  \includegraphics[width=0.7\textwidth]{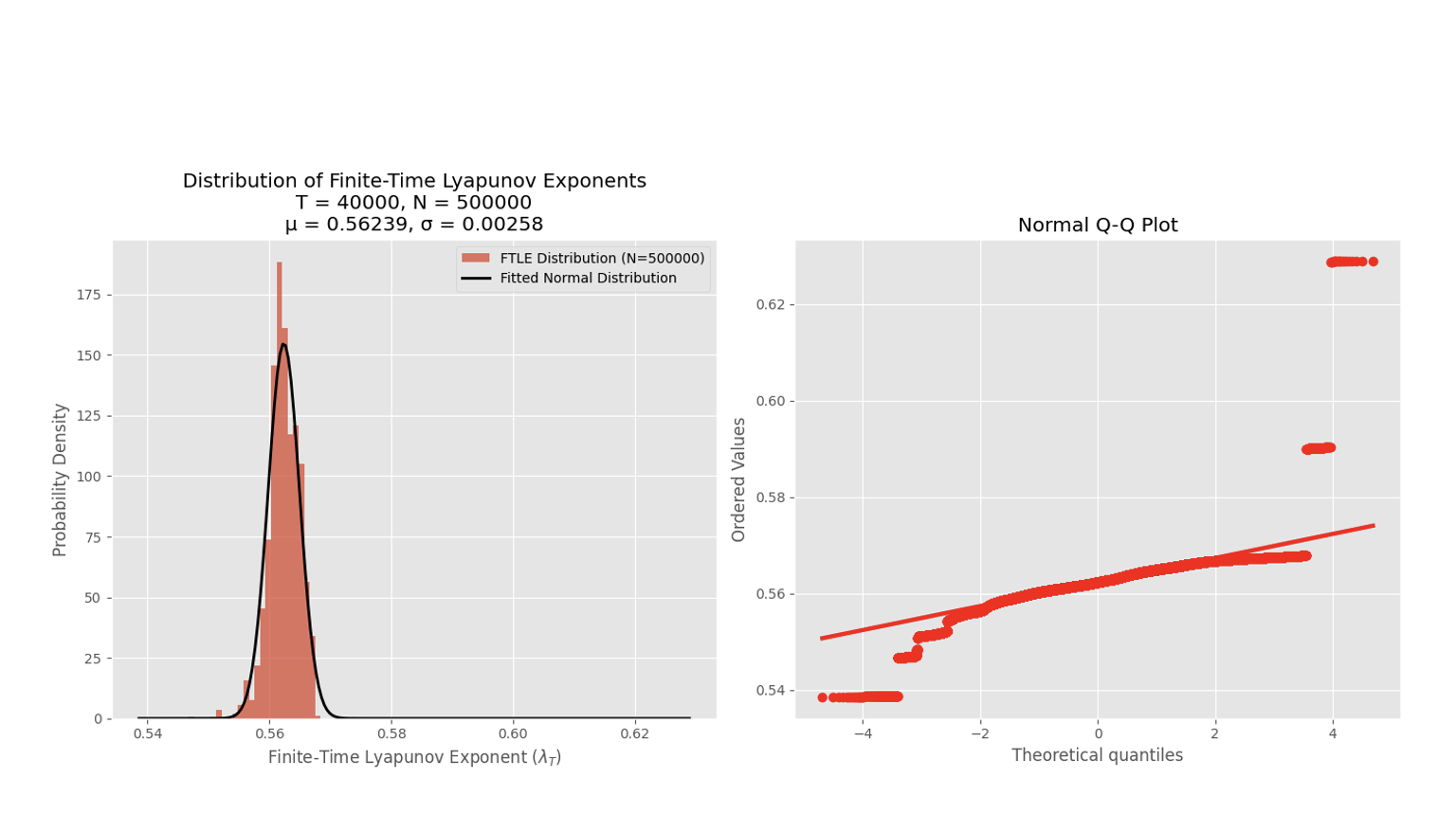}
  \caption{Histogram and Q--Q plot of the normalized FTLEs against a standard Gaussian, illustrating the excellent agreement in the bulk of the distribution.}
  \label{fig:ftle-hist}
\end{figure}

We simulated \(N = 5\times10^5\) trajectories of the Domenicali map, each of length \(T = 40{,}000\), and estimated
\(\hat\mu = 0.56239\) and \(\hat\sigma = 0.00258\).
The asymptotic Lyapunov exponent
\(\lambda_\infty \approx 0.56186\) reported in ~\cite{ElhadjSprott2008} agrees closely with \(\hat\mu\).
A Shapiro--Wilk test (Appendix~\ref{app:shapirowilk}) returns \(W = 0.9971\) with \(p \approx 3.2\times10^{-131}\), formally rejecting exact normality, but with half a million samples even minute tail deviations become statistically significant.
Crucially, the histogram and Q--Q plot (Figure~\ref{fig:ftle-hist}) display an almost perfect straight line over the central \(95\%\) of the distribution, confirming that the bulk of the finite--time Lyapunov exponents is well modeled by a Gaussian law and justifying our use of Gaussian confidence intervals in later sections.

\subsection{Impact of Noise on the Lyapunov Spectrum}

We concentrate on the two non-zero Lyapunov exponents of the
two--dimensional Domenicali map,
\[
  \lambda_1(T) \;>\; \lambda_2(T),
\qquad
  T = 20\,000 \text{ steps}.
\]
Algorithm~\ref{ftle} summarises the procedure carried out for every noise
level~$\sigma$.

\begin{algorithm}[H]
\label{ftle}
\caption{Bias--corrected finite-time Lyapunov exponents and \(\;95\%\) CIs for a single noise level~\(\sigma\)}
\label{alg:ftle-selfcontained}
\DontPrintSemicolon
\KwIn{%
  number of trajectories \(N\);              \\
  integration horizon \(T\);                 \\
  bias-correction cut--off \(T_B\) (\(1{,}000\)); \\
  noise amplitude \(\sigma\);                \\
  pre-computed \textbf{anchor clouds} \(\mathcal C_c,\mathcal C_o\) (\(\approx 300\) points each) that represent the central and outer deterministic attractors%
}
\KwOut{%
  Sample means \(\mu_{1,c},\mu_{1,o},\mu_{2,c},\mu_{2,o}\) and half--widths \(h_{1,c},h_{1,o},h_{2,c},h_{2,o}\) of the \(95\%\) confidence intervals for \(\lambda_1,\lambda_2\) in the two basins%
}

\textbf{Initialise}\;
\begin{enumerate}[leftmargin=2em,itemsep=0pt]
  \item Draw \(N\) initial states \((x_j,y_j)\) uniformly in the bounding box \([-4.5,4.5]^2\).
  \item Set tangent vectors \(v_j\gets(1,0)\) for all \(j\).
  \item Zero accumulators \(\texttt{S\_log}_j,\texttt{S\_det}_j,\texttt{S\_bias}_j\;\;(j=1,\dots,N)\).
\end{enumerate}

\textbf{Main loop}\;
\For(\tcp*[f]{iterate the noisy map}){$t=0$ \KwTo $T-1$}{
  Draw i.i.d.\ noise \(\delta a_j,\delta b_j \sim\mathcal U[-\sigma,\sigma]\) for every trajectory \(j\);\;
  Update states \((x_j,y_j)\gets f_{a+\delta a_j,b+\delta b_j}(x_j,y_j)\);\;
  Update tangent: \(w_j\gets Df_{a+\delta a_j,b+\delta b_j}(x_j,y_j)\,v_j\);\;
  \(\ell_j\gets\log\|w_j\|\);   \texttt{S\_log\(_j\)}\(+\!=\) \(\ell_j\);\;
  Normalise \(v_j\gets w_j/\|w_j\|\);\;
  \(\texttt{S\_det}_j+\!=\log|\,\det Df_{a+\delta a_j,b+\delta b_j}(x_j,y_j)|\);\;
  \If(\tcp*[f]{store bias part}){$t<T_B$}{\texttt{S\_bias\(_j\)}\(+\!=\ell_j\);   \If{$t=T_B-1$}{store snapshot \((\hat x_j,\hat y_j)=(x_j,y_j)\);}}
}

\textbf{Classify trajectories}\;
For each \(j\) compute \(d_c=\min_{q\in\mathcal C_c}\|\!(\hat x_j,\hat y_j)-q\|\) and
\(d_o\) analogously.  
Set a label \emph{central} if \(d_c<d_o\) and \(d_c<\varepsilon\) (threshold),
\emph{outer} if \(d_o<d_c\) and \(d_o<\varepsilon\).

\textbf{Compute exponents}\;
\[
\lambda_{1,j}=\frac{\texttt{S\_log}_j-\texttt{S\_bias}_j}{T},
\qquad
\lambda_{2,j}=\frac{\texttt{S\_det}_j}{T}-\lambda_{1,j}.
\]

\textbf{Aggregate statistics}\;
for each basin \(b\in\{\text{central},\text{outer}\}\) and each index \(k\in\{1,2\}\):

\begin{enumerate}[leftmargin=2em,itemsep=0pt]
  \item Collect the sample \(\{\lambda_{k,j}\}_{j\in b}\) of size \(n_b\).
  \item Compute mean \(\mu_{k,b}\) and unbiased standard deviation \(s_{k,b}\).
  \item Half--width of the \(95\%\) CI: \(h_{k,b}=1.96\,s_{k,b}/\sqrt{n_b}\).
\end{enumerate}

\Return all \(\mu_{k,b}\) and \(h_{k,b}\).
\end{algorithm}

\begin{figure}[h]
  \centering
  \begin{subfigure}[t]{0.48\textwidth}
      \centering
      \includegraphics[width=\linewidth]{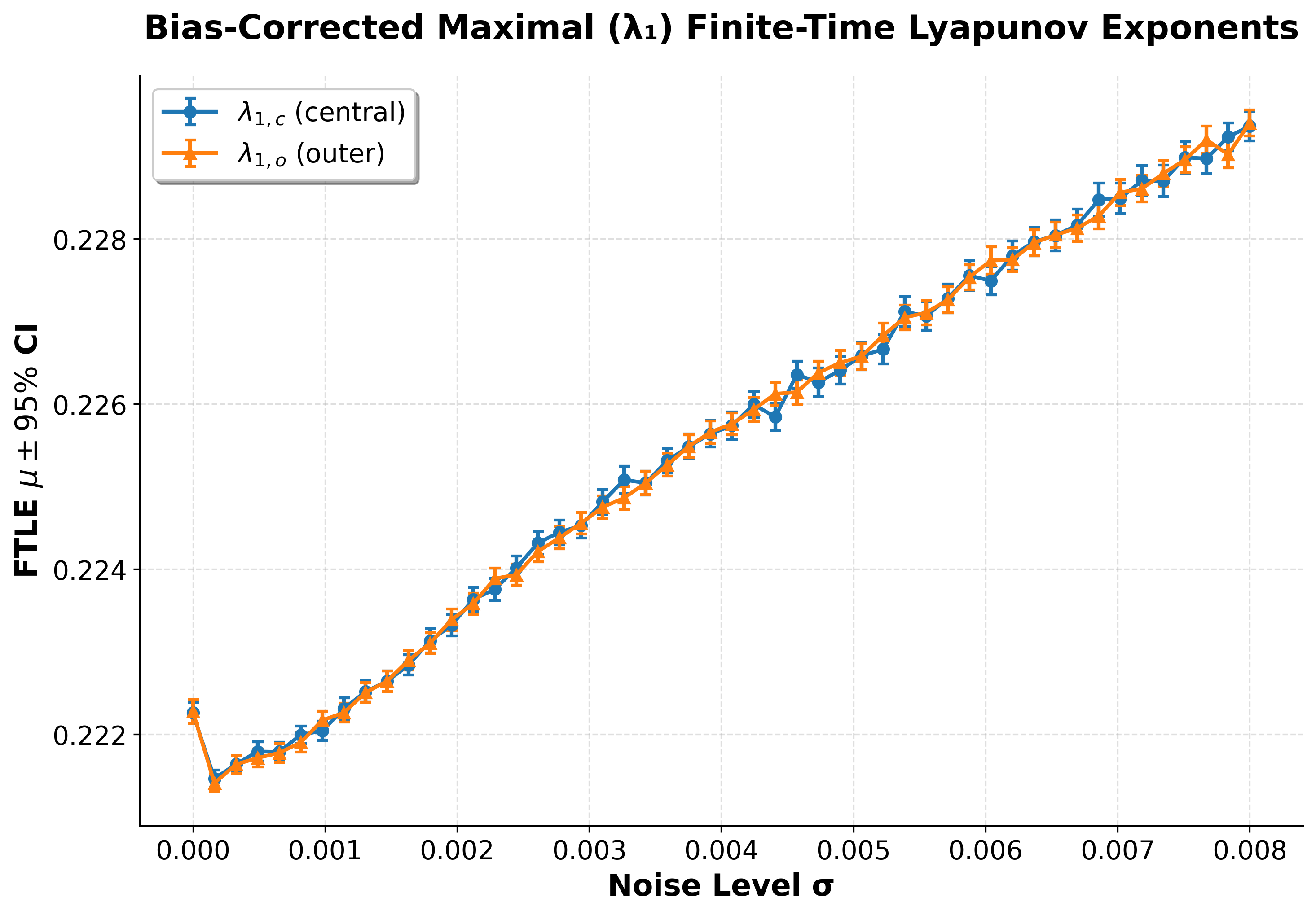}
      \caption{Maximal exponent $\lambda_1$ vs noise.}
  \end{subfigure}
  \hfill
  \begin{subfigure}[t]{0.48\textwidth}
      \centering
      \includegraphics[width=\linewidth]{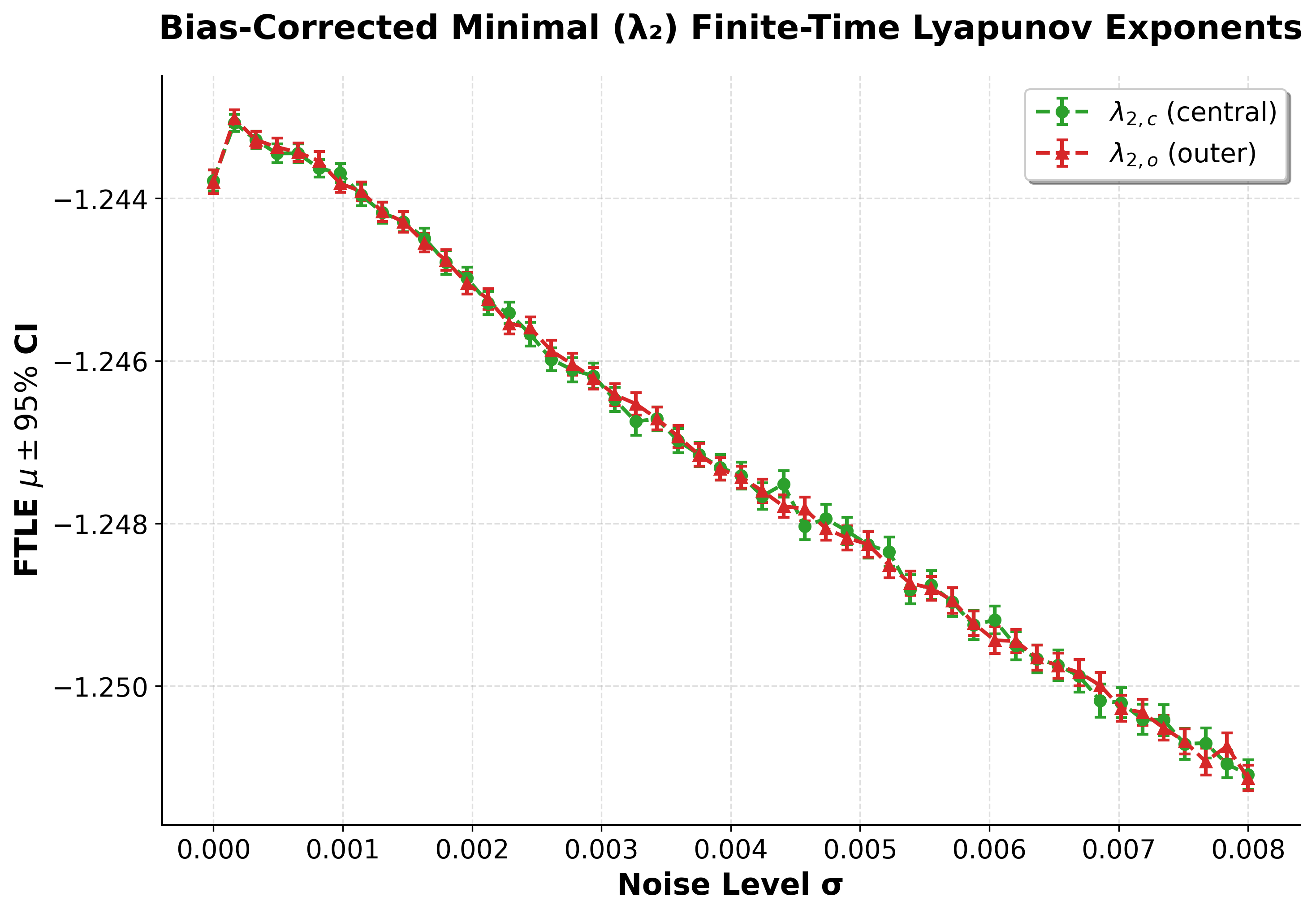}
      \caption{Minimal exponent $\lambda_2$ vs noise.}
  \end{subfigure}
  \caption{Finite-time Lyapunov exponents (mean $\pm95\%$ CI) as a function of noise amplitude $\sigma$, for both central and outer trajectory classes.}
  \label{fig:lyap_sidebyside}
\end{figure}

Figure~\ref{fig:lyap_sidebyside} presents the outcome.
Both panels share the same horizontal scale; shaded bands denote
$95\%$ confidence intervals.

The tiny dip in $\lambda_1$ and mirror bump in $\lambda_2$
for $\sigma \lesssim 2\times10^{-4}$ arise from the negative curvature
of $\lambda_1(a,b)$ at $(2.2,-0.36)$, as explained in Section~\ref{subsec:small_sigma_dip}.
Beyond that regime the two exponents diverge monotonically, widening the
spectral gap.

\vspace{1em}


We investigate how the maximal and minimal finite-time Lyapunov exponents respond to parameter noise. For each noise level $\sigma$, Algorithm~\ref{alg:sigma_run} is executed to compute $95\%$ confidence intervals (CIs) for both the central and outer attractors. We distinguish between the two attractors using the methodology in \ref{subsec:algB}. \\
\begin{algorithm}[H]
\caption{Simulation and analysis for a single noise level~$\sigma$}
\label{alg:sigma_run}
\DontPrintSemicolon
\KwIn{Noise amplitude $\sigma$, trajectories $N$, horizon $T$, burn-in $T_B$}
\KwOut{Statistics for $\lambda_1, \lambda_2$ (central/outer)}
Initialise $(x_j,y_j)$; initialise tangent $v_j=(1,0)$ and accumulators\;
\For{$t=1$ \KwTo $100$}{\tcp*{density warm-up}
  draw $\delta a,\delta b\sim\mathcal U[-\sigma,\sigma]$; evolve $(x,y)$\;
}
\For{$t=0$ \KwTo $T-1$}{\tcp*{main loop}
  sample $\delta a,\delta b$; evolve state and tangent variables\;
  accumulate $\text{sum\_log}$, $\text{sum\_det}$, $\text{bias\_log}$\;
}
classify trajectories by anchor clouds at $t=T_B$\;
\(
  \lambda_{1,j}\!=\!(\text{sum\_log}_j-\text{bias\_log}_j)/T,\;
  \lambda_{2,j}\!=\!(\text{sum\_det}_j/T)-\lambda_{1,j}
\)\;
compute means $\mu_{1,2}$ and CIs
$\mu_i\pm z_{.975}s_i/\sqrt{n_i}$ for both groups\;
\end{algorithm}

\subsubsection{Small-$\sigma$ dip and eventual divergence}
\label{subsec:small_sigma_dip}

Consider the deterministic spectrum as a smooth map
\(F:(a,b)\mapsto(\lambda_{1},\lambda_{2})\).
At the base point \((a_{0},b_{0})=(2.2,-0.36)\) we perturb the parameters with zero-mean i.i.d.\ noise \((\delta a,\delta b)\sim\mathcal U[-\sigma,\sigma]^2\).
A second-order Taylor expansion gives
\[
\lambda_{1}(a_{0}+\delta a,b_{0}+\delta b)
   =\lambda_{1}^{0}
    +\nabla\lambda_{1}^{0}\!\cdot\!(\delta a,\delta b)
    +\tfrac12(\delta a,\delta b)H_{1}^{0}(\delta a,\delta b)^{\!\top}
    +\mathcal O(\|\delta\|^{3}),
\]
where \(H_{1}^{0}\) is the Hessian at the base point [\ref{app:taylor}].
Taking expectations,
\[
\mathbb E\bigl[\lambda_{1}(\sigma)\bigr]
  =\lambda_{1}^{0}
   +\tfrac12\operatorname{tr}\!\bigl(H_{1}^{0}\,\mathbb E[\delta\delta^{\top}]\bigr)
   +\mathcal O(\sigma^{3})
  =\lambda_{1}^{0}
   +\frac{\sigma^{2}}{6}\operatorname{tr}(H_{1}^{0})
   +\mathcal O(\sigma^{3}),
\]
because \(\mathbb E[\delta a]=\mathbb E[\delta b]=0\) and
\(\mathbb E[\delta^{2}]=\sigma^{2}/3\).
At \((2.2,-0.36)\) the trace is negative, hence
\(\Delta\lambda_{1}<0\):
Jensen's inequality [\ref{app:jensen}] forces the ensemble-averaged \(\lambda_{1}\) to dip for small noise.

Beyond the perturbative window (\(\sigma\gtrsim10^{-3}\)) higher-order terms dominate; trajectories begin to sample regions where \(\lambda_{1}\) is intrinsically larger, so \(\lambda_{1}\) rises and \(|\lambda_{2}|\) deepens, widening the spectral gap by \(\approx0.007\) over \(\sigma\in[0,8\times10^{-3}]\).

By similar arguments, the same applies for \(\lambda_{2}\): the opposite curvature of its surface produces the mirror ``bump'' at small \(\sigma\), followed by a monotone decrease as the noise grows.

\subsubsection{Noise-Induced Homogenisation}
By the distribution of the FTLEs, we have the variance shrinking like $\frac{1}{T}$, thus by $T = 20\,000$ the spread is already below our sampling error. In other words, the bias from burn-in becomes negligible at large $T$, implying that the residual overlap of the 95\% confidence intervals in Figure~\ref{fig:lyap_sidebyside} is a genuine dynamical homogenisation (all trajectories sample essentially the same mix of expanding and contracting stretches), not a consequence of inadequate statistics.

\subsubsection{Kaplan--Yorke Dimension and Error Propagation}

For two dimensions,
\(
  D_{\mathrm{KY}} = 1 + \frac{\lambda_1}{|\lambda_2|}.
\)
Treating $(\lambda_1,\lambda_2)$ as jointly normal, we propagate uncertainty with the multivariate delta method [\ref{app:mdelta}]. Writing $g(\lambda_1,\lambda_2)=1+\lambda_1/|\lambda_2|$,
\[
  \partial_{\lambda_1}g=\frac{1}{|\mu_2|},\qquad
  \partial_{\lambda_2}g=-\frac{\mu_1}{\mu_2^{2}},
\]
\[
  \operatorname{Var}(D_{\mathrm{KY}}) \approx
    \frac{\operatorname{Var}(\lambda_1)}{\mu_2^{2}}
  + \frac{\mu_1^{2}}{\mu_2^{4}}\operatorname{Var}(\lambda_2)
  - 2\frac{\mu_1}{\mu_2^{3}}\operatorname{Cov}(\lambda_1,\lambda_2),
  \qquad
  \text{CI}_{95\%}: \;
  D_{\mathrm{KY}} \pm z_{0.975}\sqrt{\operatorname{Var}(D_{\mathrm{KY}})}.
\]
\begin{figure}[H]
  \centering
  \includegraphics[width=.7\linewidth]{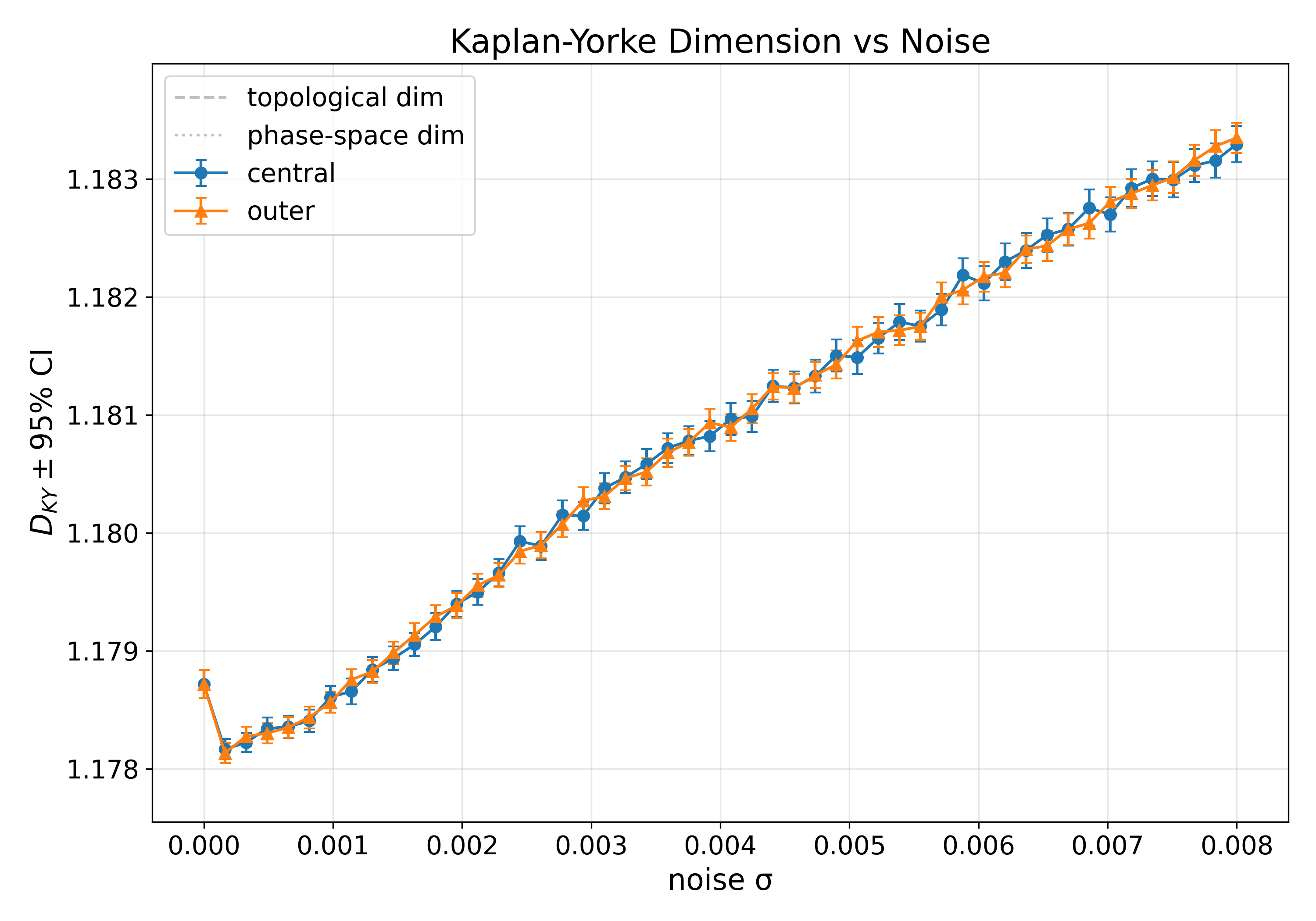}
  \caption{Kaplan--Yorke dimension with $95\%$ CIs.
           The ordinate is zoomed to reveal the statistically
           significant yet modest $\approx0.004$ increase over the noise range.}
  \label{fig:kydim_noise}
\end{figure}
\FloatBarrier

The gentle upward drift of \(D_{\mathrm{KY}}\) reveals how parameter noise broadens the invariant set without changing its topological class.  In the deterministic limit the attractor behaves like a thin filament with \(D_{\mathrm{KY}}\simeq1.178\).  Adding noise thickens this filament into a fuzzy tube: the leading exponent grows, the contracting one becomes more negative, and their ratio \( \lambda_{1}/|\lambda_{2}| \) increases.  Over the explored range (\(0\le\sigma\le0.008\)) the dimension rises by \(\approx0.004\), corresponding to an extra \(4\times10^{-3}\) effective degrees of freedom.  The growth is nearly linear in \(\sigma\) and, more interestingly, the $D_{KY}$ increase was the same whether we looked at trajectories in either attractor. This reinforces that noise is causing a uniform thickening of the attractors' combined support in phase space.

\newpage

\section{Numerical Investigation of Noise Effects}
\label{sec:numerics} 

In this section we describe three complementary numerical algorithms developed to study how additive parametric noise influences the structure and stability of coexisting attractors in the Domenicali map.  First, we present a grid-based histogram method that yields high-resolution density fields.  Second, we present a method to produce accurate deterministic basins of attraction efficiently, using a reference-cloud of anchor points. Finally, we introduce a two-pass approach to produce accurate random basins, making use of a deterministic basin from the second algorithm, augmented by a stochastic ``random basin'' pass and a final reference-cloud, batched verification step. 

For quick development, we opted for \texttt{Python}, and used \texttt{Google} \texttt{Colab} for free and powerful GPU compute. This combination, together with the above tools allow us to locate basin boundaries, quantify noise-induced ``collisions'' between random attractors, and estimate a critical noise upper bound \(\sigma_c\le 4.5\times10^{-3}\) at which a collision occurs, both efficiently and in a robust manner.

\subsection{Algorithm A: High-Resolution Histogram Density Field}
\label{subsec:algA}

\subsubsection{Purpose}  Visualise the geometry of the noisy attractor as a smooth density field.

\subsubsection{Inputs and Parameters}
\begin{itemize}
  \item Number of sample points: \(N_{\text{pts}} = 50\,000\).
  \item Grid resolution: \(n_x = n_y = 800\).
  \item Window: \([-4.5,4.5]\times[-4.5,4.5]\).
  \item Noise amplitude \(\sigma\in[0,0.008]\).
  \item Time steps per trajectory: \(T = 50\,000\).
  \item Burn-in: \(T_{\mathrm{burn}} = 1\,000\).
\end{itemize}

\subsubsection{Procedure}
\begin{enumerate}
  \item Initialise \(N_{\text{pts}}\) points uniformly in the window.
  \item Evolve each point through \(T_{\mathrm{burn}}\) steps under
    \[
      \begin{cases}
        a_t = 2.2 + \delta_{a,t},\quad \delta_{a,t}\sim \mathcal{U}[-\sigma,\sigma],\\
        b_t = -0.36 + \delta_{b,t},\quad \delta_{b,t}\sim \mathcal{U}[-\sigma,\sigma],
      \end{cases}
    \]
    to ensure convergence into the attractors.
  \item Zero an \(n_x\times n_y\) integer array \(\mathcal{H}\).
  \item For \(t = 1,\dots,T\):
    \begin{enumerate}
      \item Evolve each point \((x,y)\mapsto(1 - a_t\sin x + b_t y,\;x)\).
      \item Map \((x,y)\) back to grid-indices \((i,j)\) by interpolation.
      \item Increment \(\mathcal{H}_{i,j}\leftarrow\mathcal{H}_{i,j}+1\).
    \end{enumerate}
  \item After \(T\) steps, normalise \(\mathcal{H}\) to obtain a density field.
  \item Render \(\mathcal{H}\) with a smooth colormap (e.g.\ \texttt{inferno}).
\end{enumerate}
\begin{figure}[H]
  \centering
  \caption*{\textbf{Evolution of Attractors Density with Increasing Noise}}
  \begin{subfigure}[b]{0.48\textwidth}
    \centering
    \includegraphics[width=\linewidth]{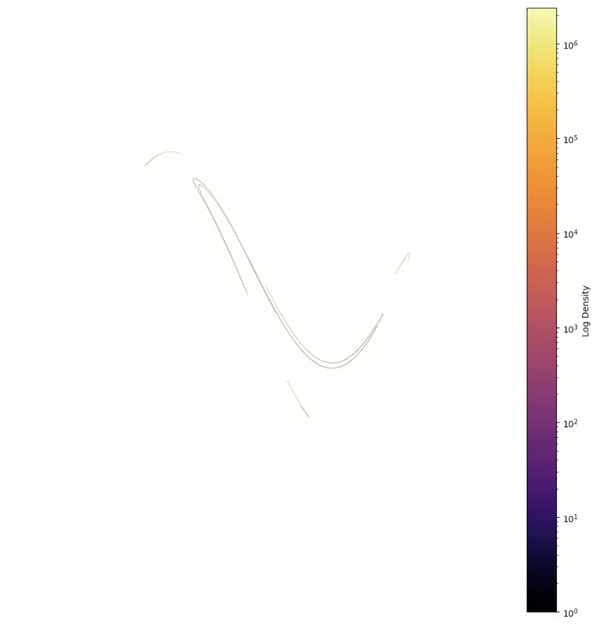}
    \caption{Small noise (\(\sigma = 0.0001\)).}
    \label{fig:density-sigma0040}
  \end{subfigure}
  \hfill
  \begin{subfigure}[b]{0.48\textwidth}
    \centering
    \includegraphics[width=\linewidth]{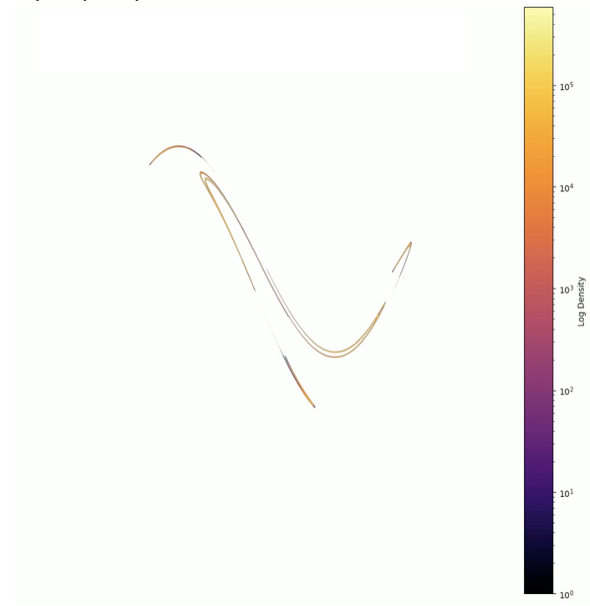}
    \caption{Moderate noise (\(\sigma = 0.0044\)).}
    \label{fig:density-sigma0500}
  \end{subfigure}

  \vspace{1em}

  \begin{subfigure}[b]{0.48\textwidth}
    \centering
    \includegraphics[width=\linewidth]{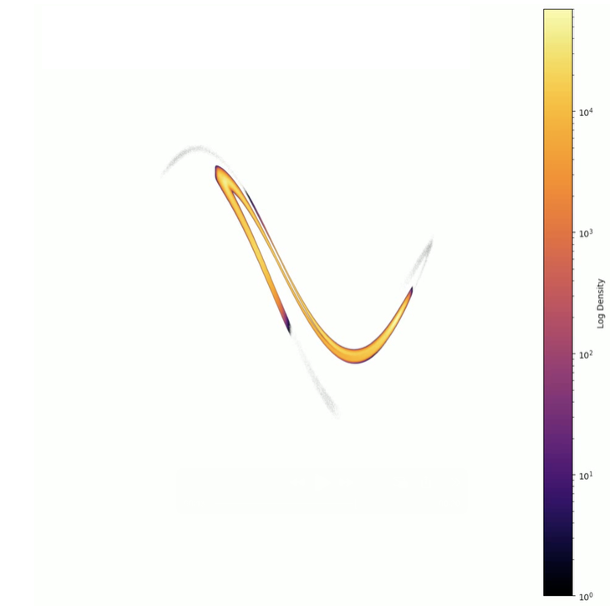}
    \caption{Large noise (\(\sigma = 0.0435\)).}
    \label{fig:density-sigma0100}
  \end{subfigure}
  \hfill
  \begin{subfigure}[b]{0.48\textwidth}
    \centering
    \includegraphics[width=\linewidth]{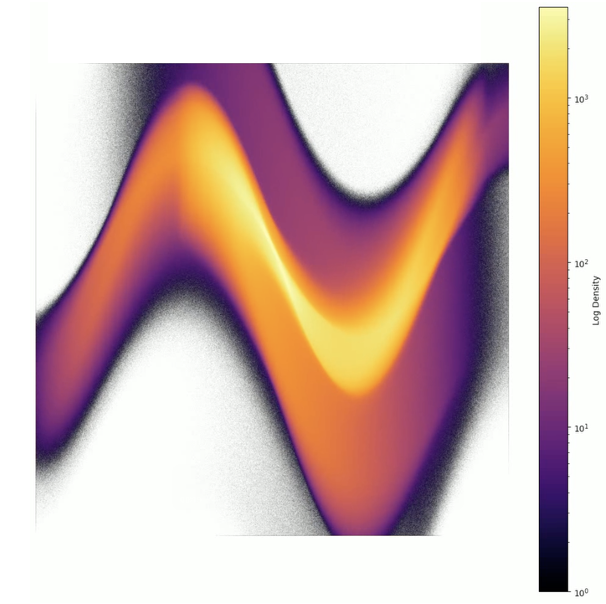}
    \caption{Very large noise (\(\sigma = 0.5706\)).}
    \label{fig:density-sigma1000}
  \end{subfigure}
\caption{High-resolution density fields of the Domenicali attractor generated with Algorithm A.}
  \label{fig:hist_density}

  \vspace{0.5em}
  \begin{enumerate}[label=(\alph*)]
    \item With almost no noise, the two filaments remain distinct.
    \item The outer attractor begins to touch the inner one.
    \item The outer attractor starts to vanish as most points converge onto the middle filament.
    \item Under very high noise, the attractor thickens substantially yet stays bounded.
  \end{enumerate}
\end{figure}

\newpage

\subsection{Algorithm B: Reference-Based Deterministic Basins}
\label{subsec:algB}
\subsubsection{Purpose}
To efficiently partition the phase space into two distinct basins of attraction for the deterministic Domenicali map (\(\sigma=0\)), we employ a \emph{reference-cloud} methodology instead of computationally intensive clustering. Although chaotic dynamics amplify small perturbations over short timescales, trajectories that converge to the same attractor will, in the long-term, recurrently approach one another arbitrarily closely. By numerically sampling two finite sets of ``anchor'' points, each drawn from a trajectory initiated in the basin of a known attractor, from a particular seed, we obtain two reference clouds. We then classify any arbitrary initial condition by computing its smallest Euclidean distance to a point in each cloud and assigning it to the basin corresponding to the anchor set it has gotten close to. 


\medskip
\subsubsection{Inputs and Parameters}
\begin{itemize}
  \item \((a,b)=(2.2,-0.36)\): deterministic map parameters.
  \item Two seed points \(\,(x,y) = (0,0)\) (central attractor) and \((1.5,1.0)\) (outer attractor).  
  \item \(m_{\rm ref}=300\): number of reference samples per cloud.  
  \item \(T_{\rm det}=50{,}000\): number of evolution steps.  
  \item \(\varepsilon=0.0044\): proximity threshold for classification.  
  \item Grid of \(g\times g\) initial points (\(g=2000\)) over \([-4.5,4.5]^2\).  
\end{itemize}

\medskip
\subsubsection{Outputs}
\begin{itemize}
  \item Index-sets \(\mathit{idx}_1,\mathit{idx}_2\) and label tensor \texttt{labels} of length \(N=g^2\), with values 
  \(\{0=\text{undefined},1=\text{central},2=\text{outer}\}\).
  \item Scatter plot showing red/black (two basins) and yellow (undefined) points.
\end{itemize}

\medskip
\subsubsection{Overview of steps}
\begin{enumerate}
  \item \textbf{Generate Anchor Clouds.}  
    \begin{enumerate}
      \item For each seed, evolve the deterministic map for a short burn-in.
      \item Continue for \(m_{\rm ref}\) steps, collecting each \((x,y)\) into a small tensor \(\mathcal{C}_c\) or \(\mathcal{C}_o\).
    \end{enumerate}
  \item \textbf{Build Phase-Space Grid.}  
    Generate \(N=g^2\) equally spaced initial conditions \(\{(x_{0,i},y_{0,i})\}\).
  \item \textbf{Parallel Evolution \& Classification.}  
    Initialise \(\texttt{labels}[:] = 0\).  Then for \(t=1\ldots T_{\rm det}\):
    \begin{itemize}
      \item Update all points \((x,y)\mapsto(1 - a\sin x + b\,y,\;x)\).
      \item Identify still-undefined indices \(\mathit{idx}_{\mathrm{und}}\).
      \item Stack their positions \(\mathbf{P}\in\mathbb{R}^{M\times2}\).
      \item Compute distances \(d_1 = \min_{q\in\mathcal{C}_c}\|\mathbf{P}-q\|\), and \(d_2\) similarly.
      \item Label as central if \(d_1<\varepsilon\le d_2\), or outer if \(d_2<\varepsilon\le d_1\).
    \end{itemize}
  \item \textbf{Extract \& Plot.}  Gather \(\mathit{idx}_1,\mathit{idx}_2\), and plot colored scatter.
\end{enumerate}

\medskip


\begin{figure}[H]
  \centering
  \includegraphics[width=.7\linewidth]{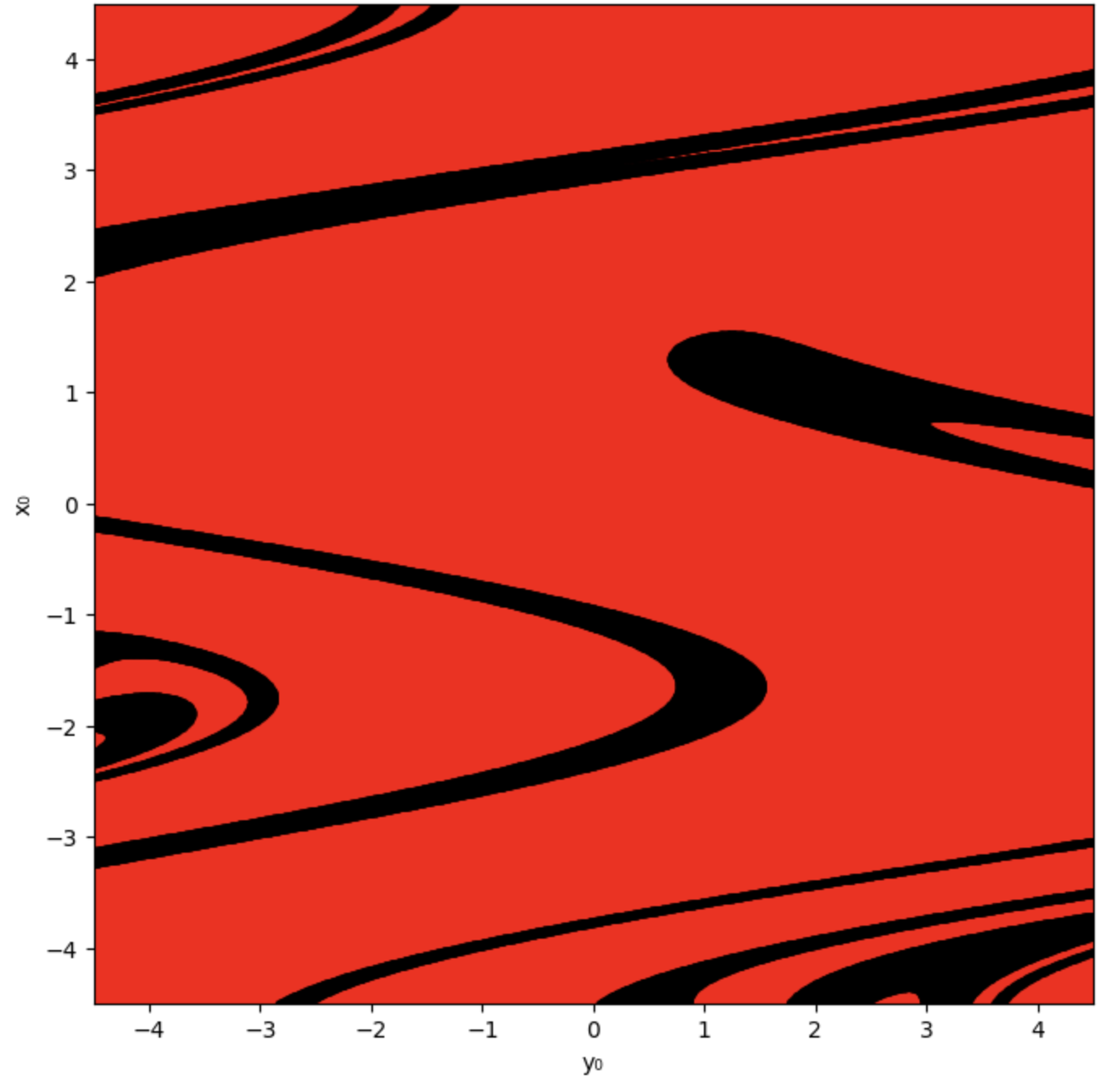}
  \caption{Scatter of initial conditions coloured by
  deterministic basin. The red region converges to the visually middle attractor, and the black region converges to the outer period-3 attractor}
  \label{fig:deterministic_basin}
\end{figure}

\medskip
\subsubsection{Efficiencies in the code}
\begin{itemize}
  \item \emph{Vectorised evolution:} All \(M\) remaining points are updated in one batched tensor operation per step.
  \item \emph{Shrinking checks:} Only undefined points are distance-checked, which quickly diminishes in our deterministic basin.
\end{itemize}

\subsection{Algorithm C: Multi-Step Random Basins}
\label{subsec:algC}

\subsubsection{Challenges in Defining and Computing Random Basins}


In the deterministic setting, the \emph{basin of attraction} of an attractor \(A\subset X\) is commonly defined as
\[
  B(A)\;=\;\{\,x_0\in X : \omega(x_0)\subset A\},
\]
where \(\omega(x_0)\) denotes the \emph{\(\omega\)-limit set} of the trajectory starting at \(x_0\).  Numerically one simply colors each grid-point by which attractor its orbit lands in.

Under \emph{parametric noise}, however, several new difficulties arise:
\begin{itemize}
  \item \textbf{Non-invariance of orbits:}  A noisy trajectory can hop between basins arbitrarily late, so it may not settle permanently into any one attractor.
  \item \textbf{Sample-path dependence:}  Each realisation of the noise yields a different map and hence a different \(\omega\)-limit set.
  \item \textbf{Ambiguity of ``eventual'' behavior:}  A point might spend long times near one attractor before switching, how should its color be chosen?
\end{itemize}

\medskip
\noindent To address this, we adopt a \emph{random-dynamical-systems} viewpoint.  Let \((\Omega,\mathcal{F},\mathbb{P})\) be the noise probability space, \(\Theta_t\) the shift on noise realisations, and
\[
  \Phi(t,\omega)\;:\;X\to X
\]
the stochastic flow generated by the noisy Domenicali map. A \emph{random attractor} is a measurable set \(A(\omega)\subset X\) satisfying
\[
  \Phi(t,\omega)\bigl(A(\omega)\bigr)\;=\;A\bigl(\Theta_t\omega\bigr),
  \quad\text{and}\quad
  \lim_{t\to\infty}\mathrm{dist}\bigl(\Phi(t,\Theta_{-t}\omega)x,\;A(\omega)\bigr)=0
  \quad\forall\,x\in B(A)(\omega).
\]
Here the \emph{random basin} of \(A\) at realisation \(\omega\) is defined by
\[
  B(A)(\omega)
  \;=\;
  \bigl\{\,x\in X : \lim_{t\to\infty}\mathrm{dist}\bigl(\Phi(t,\Theta_{-t}\omega)x,\;A(\omega)\bigr)=0 \bigr\}.
\]
Equivalently, one can use the \(\omega\)-limit set in the skew-product:
\[
  \omega_\omega(x)
  \;=\;
  \bigcap_{T\ge0}\overline{\bigl\{\Phi(t,\Theta_{-t}\omega)x : t\ge T\bigr\}},
  \quad
  B(A)(\omega)
  =\;
  \{\,x : \omega_\omega(x)\subset A(\omega)\}.
\]
This formalism captures ``eventual almost-sure convergence'' under the same noise realisation.

\medskip
\subsubsection{Numerical implications}
\begin{enumerate}
  \item One must decide on a \emph{finite-time approximation} of \(\omega_\omega(x)\), e.g.\ sample the last \(T\) iterates.
  \item Points whose sample-path \(\omega\)-limits switch between attractors require a special ``white'' label.
  \item Efficient computation again exploits \emph{reference clouds} for each random attractor and on-the-fly distance checks (see Algorithms B [\ref{subsec:algB}] above and C [\ref{subsec:algC}] below).
\end{enumerate}

\subsubsection{Overview}
Algorithm C constructs a deterministic basin (refer to Algorithm B [\ref{subsec:algB}]), and then processes the set of points in the deterministic basin by evolving them through the random Domenicali map, checking at each time step if the evolved points have stayed, or left their parent deterministic basins. Lastly, a final check to ensure "white points" (those points that have been observed to switch deterministic basins under the random mapping) have been classified correctly, by evolving them and seeing if they ever get close to a reference cloud of a known trajectory.

\subsubsection{Step 1: Reference Deterministic Basin}
Refer to Algorithm B [\ref{subsec:algB}].

\subsubsection{Step 2: Random-Basin Noisy Pass}
\paragraph{Inputs}
\begin{itemize}
  \item Noise amplitude \(\sigma\in(0,0.008]\).  
  \item Deterministic basin index-sets \(\idx_{\rm central}\) and \(\idx_{\rm outer}\) (each a subset of \(\{1,\dots,N\}\)).  
  \item Total random-evolution length \(T_{\rm rnd}=10\,000\).  
    \item Grid size for our evolved points. We set to 2000 to align with Algorithm B [\ref{subsec:algB}], giving number of points \(N = {2000}^2\). 
  \item Grid-to-label lookup: a precomputed tensor \(\mathtt{basin\_map}[1\!:\!N]\in\{1,2\}^N\) giving the deterministic label of each of the \(N\) grid-points, passed during Step 1.
\end{itemize}

\paragraph{State Variables}
\begin{itemize}
  \item \(\mathbf{x},\mathbf{y}\in\mathbb{R}^N\): current coordinates of the \(N\) points (initially at grid positions).  
  \item \(\mathit{alive}\in\{0,1\}^N\): a Boolean mask (tensor) with \(\mathit{alive}_i=1\) iff point \(i\) has not yet crossed.  
  \item \(\mathit{crossed}\in\{0,1\}^N\): a Boolean mask with \(\mathit{crossed}_i=1\) once point \(i\) leaves its original basin.  
  \item \(\ell_0\in\{1,2\}^N\): the original deterministic label for each point, copied from \(\mathtt{basin\_map}\).
\end{itemize}

\paragraph{Procedure}
\begin{enumerate}
  \item \textbf{Initialization:}
    \begin{itemize}
      \item For each index \(i\in\idx_{\rm central}\cup\idx_{\rm outer}\), set \((\mathbf{x}_i,\mathbf{y}_i)\leftarrow(x_{0,i},y_{0,i})\).  
      \item Set \(\mathit{alive}_i\leftarrow1\) and \(\mathit{crossed}_i\leftarrow0\) for all \(i=1,\dots,N\).  
      \item Store original labels \(\ell_{0,i} \leftarrow \mathtt{basin\_map}_i\).
    \end{itemize}

  \item \textbf{Time-stepping:}  For \(t=1,2,\dots,T_{\rm rnd}\):
    \begin{enumerate}
      \item \emph{Evolve the alive points:}  
        \[
          a_t = a + \delta_{a},\quad b_t = b + \delta_{b},\quad
          x_{i}\mapsto 1 - a_t\sin(x_i) + b_t\,y_i,\quad y_{i}\mapsto x_i
        \]
        only for those \(i\) with \(\mathit{alive}_i=1\).
      \item \emph{Lookup current label:}  
        \[
          \ell_{t,i} \;=\; \mathtt{basin\_map}_j
          \quad\text{where }j=\mathrm{round}\bigl((x_i,y_i)\mapsto\text{grid index}\bigr).
        \]
      \item \emph{Detect crossings:}  
        Find all indices
        \(\{\,i: \mathit{alive}_i=1 \text{ and } \ell_{t,i}\neq \ell_{0,i}\}\).  For each such \(i\):
        \[
          \mathit{crossed}_i \leftarrow 1,\quad \mathit{alive}_i\leftarrow 0.
        \]
      \item \emph{Early exit:}  If \(\sum_i \mathit{alive}_i = 0\), break the loop.
    \end{enumerate}
\end{enumerate}

\begin{figure}[H]
  \centering
  \begin{subfigure}[b]{0.48\textwidth}
    \centering
    \includegraphics[width=\linewidth]{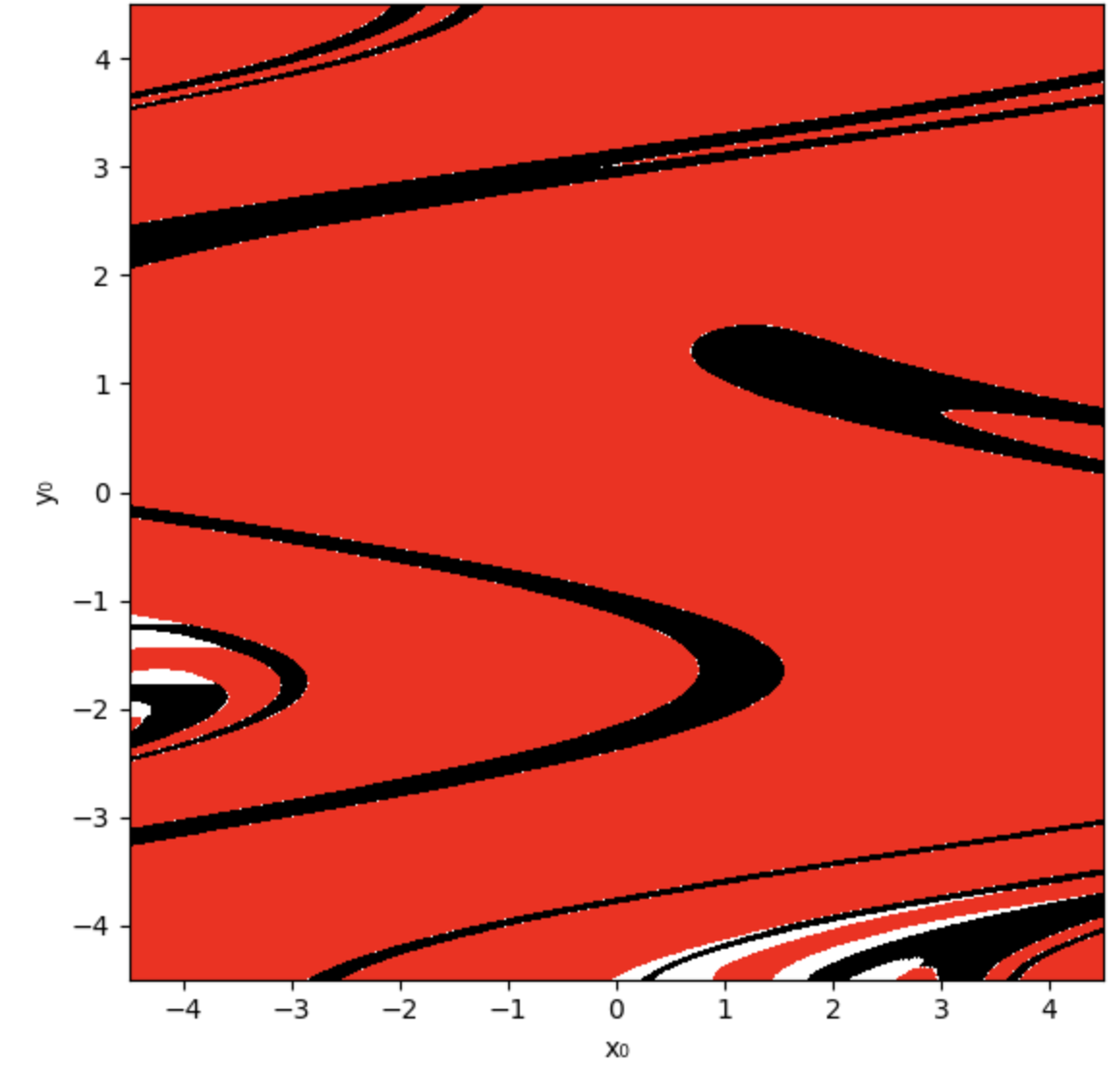}
    \caption{Small noise (\(\sigma = 5\times10^{-4}\)).}
    \label{fig:basin-single-sigma-0.0005}
  \end{subfigure}
  \hfill
  \begin{subfigure}[b]{0.48\textwidth}
    \centering
    \includegraphics[width=\linewidth]{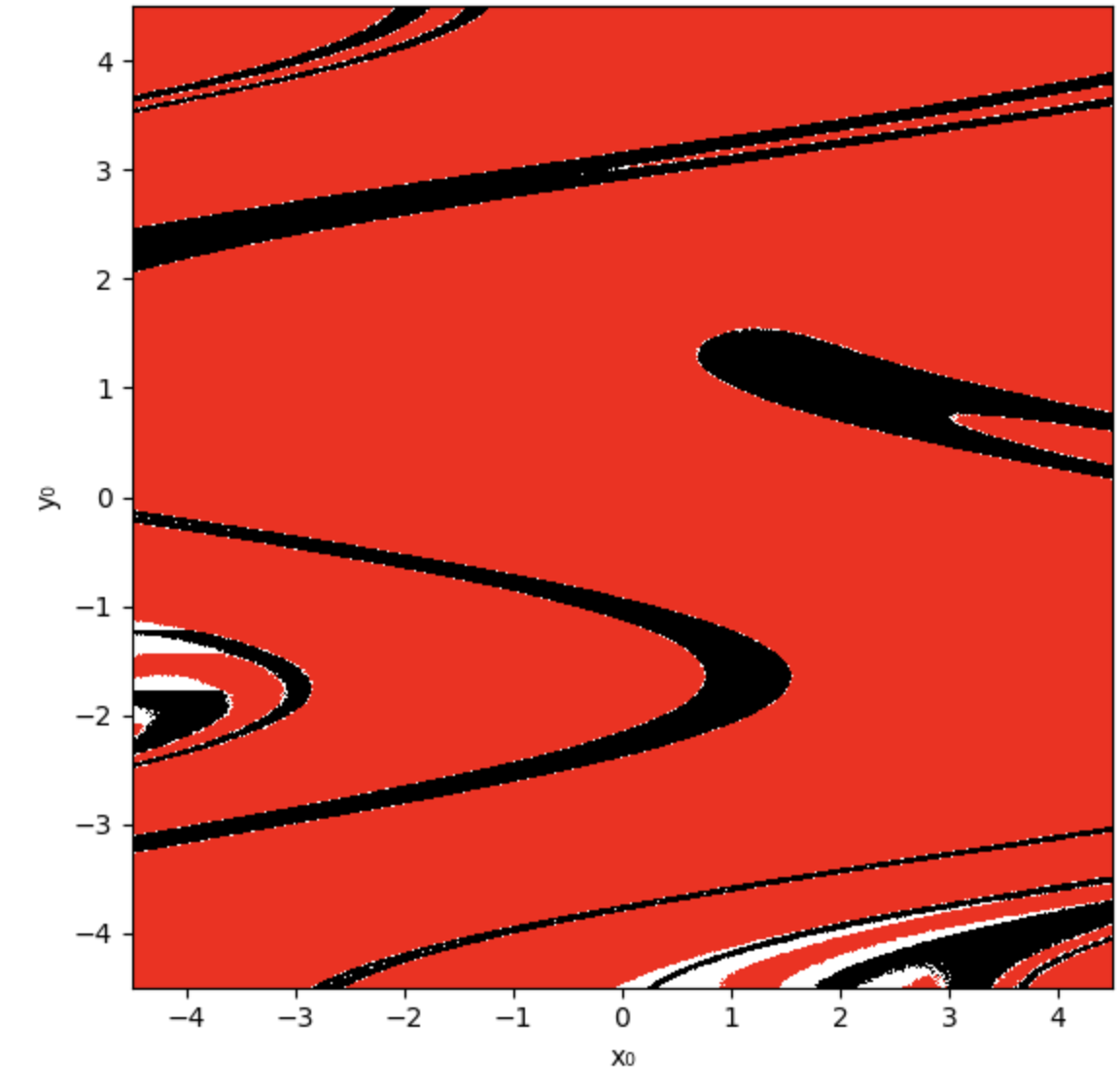}
    \caption{Moderate noise (\(\sigma = 2\times10^{-3}\)).}
    \label{fig:basin-single-sigma-0.0020}
  \end{subfigure}

  \caption{Random basin plots for the Domenicali map. A thin, fuzzy and white border can be seen at the boundary of the two basins, which gets "fatter" with more noise. There also appears unusually large white patches in some areas, possibly because of rounding errors, and exceptional sensitivity in those areas.}
  \label{fig:random_basin_initial}
\end{figure}

\subsubsection{Step 3: Dynamic Re-check}
\textbf{Purpose}\\
Remove \emph{false positive} basin-switching points (due to rounding/calculation errors).

\paragraph{Inputs}
\begin{itemize}
  \item Boolean mask \(\mathit{crossed}\in\{0,1\}^N\) from Step 2.  
  \item Second evolution step length \(T_{\rm dyn}=2{,}000\).  
  \item Same proximity threshold \(\varepsilon\) and reference clouds \(\C_c,\C_o\) from Step 1.
\end{itemize}

\paragraph{Procedure}
\begin{enumerate}
  \item Gather the list of indices
    \(\idx_{\rm crossed}=\{\,i:\mathit{crossed}_i=1\}\).
  \item Generate a single noise-sequence \(\{\delta_t\}_{t=1}^{T_{\rm dyn}}\) once.
  \item Partition \(\idx_{\rm crossed}\) into batches of size \(B\):
    \begin{itemize}
      \item For each batch, initialise local coordinate arrays \(\mathbf{x}_b,\mathbf{y}_b\).  
      \item Set a local Boolean mask \(\mathit{confirmed}\in\{0,1\}^B\) to zero.
      \item For \(t=1\ldots T_{\rm dyn}\):
        \begin{itemize}
          \item Evolve all \(B\) points under the \emph{same} \(\delta_t\).  
          \item For each attractor cloud \(k\in\{c,o\}\), compute distances
            \(\min_{q\in \C_k}\|p_i - q\|\) for as yet un\-confirmed points.
          \item If \(\min\|\cdot\|<\varepsilon\), mark that point \emph{confirmed} for cloud \(k\).  
          \item Break early if all \(B\) are confirmed.
        \end{itemize}
      \item Any batch-index \(j\) with \(\mathit{confirmed}_j=0\) was a false positive:  
        set \(\mathit{crossed}_{i_j}\leftarrow0\).
    \end{itemize}
  \item The remaining indices with \(\mathit{crossed}_i=1\) are the true noise-induced basin transitions.
\end{enumerate}

\begin{figure}[H]
  \centering
  \begin{subfigure}[b]{0.48\textwidth}
    \centering
    \includegraphics[width=\linewidth]{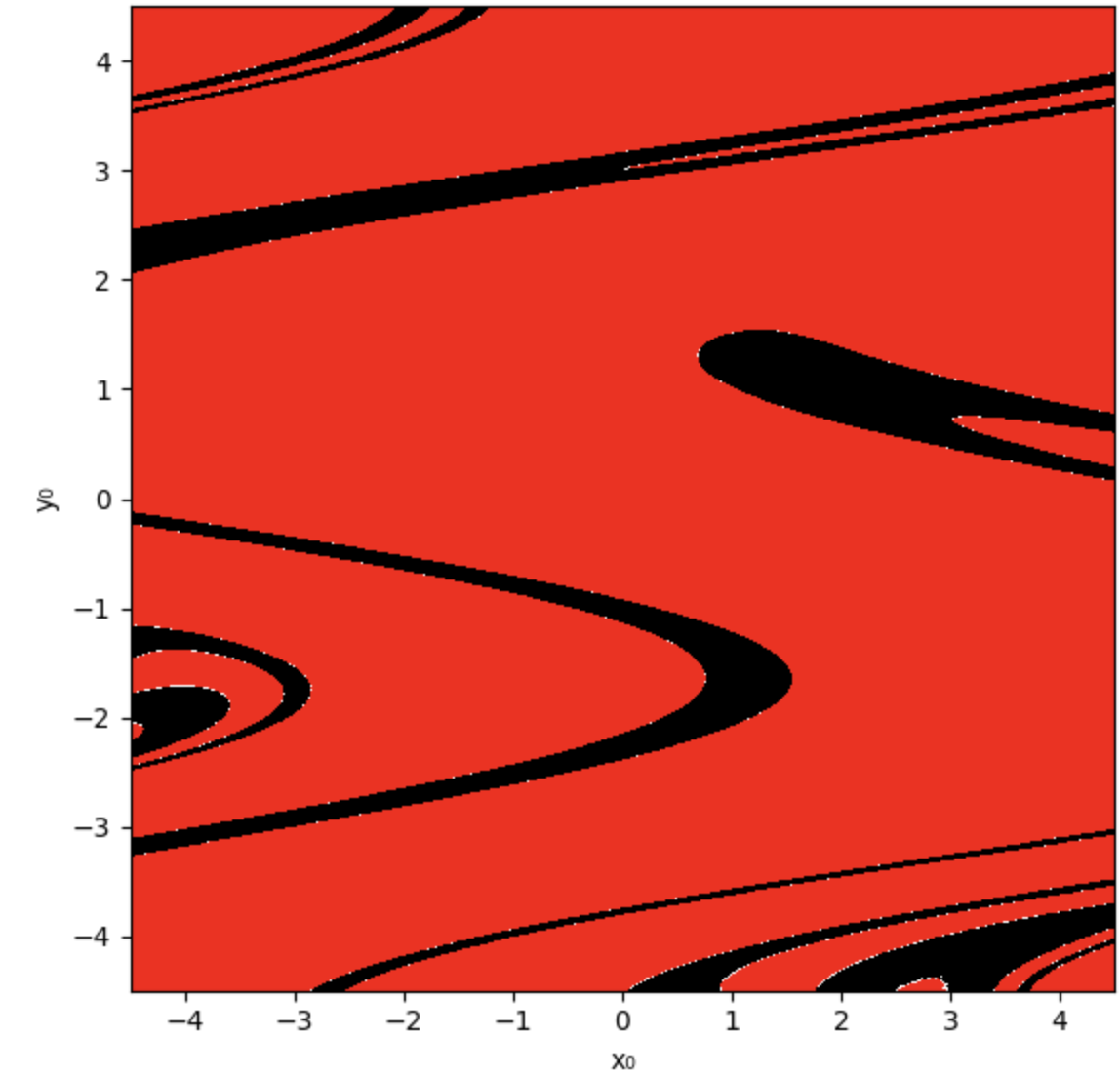}
    \caption{Small noise (\(\sigma = 5\times10^{-4}\)).}
    \label{fig:basin-double-sigma-0.0005}
  \end{subfigure}
  \hfill
  \begin{subfigure}[b]{0.48\textwidth}
    \centering
    \includegraphics[width=\linewidth]{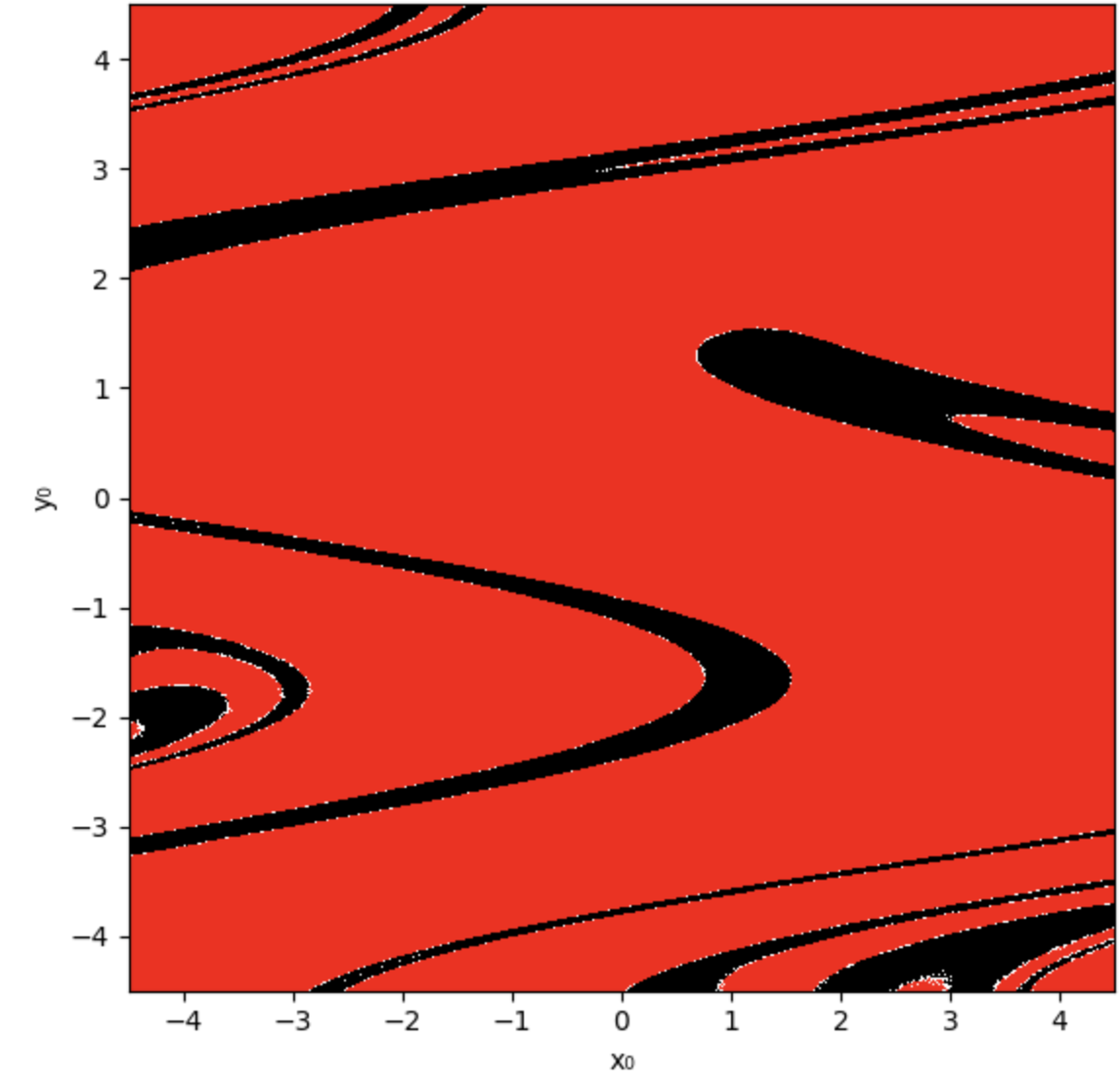}
    \caption{Moderate noise (\(\sigma = 2\times10^{-3}\)).}
    \label{fig:basin-double-sigma-0.0020}
  \end{subfigure}

  \caption{Random basin plots for the Domenicali map. This time, a secondary pass verifies the white points, and discards those that are false. The white border is still visible and still gets "fatter". The white patches from before are no longer visible either.}
  \label{fig:random_basin_verified}
\end{figure}

\paragraph{Note:} This algorithm does not evolve points for many realisations. Thence, it is most limited by both the number of realisations, as well as the number of iterations in both steps 2 and 3.

\subsection{Complexity Comparisons of Algorithms}

\begin{table}[h]
\centering
\begin{tabular}{|l|l|l|}
\hline
Algorithm & Best--case time & Worst--case time \\ \hline
A -- Histogram density & \multicolumn{2}{c|}{$O\!\bigl(N\,T\bigr)$} \\ \hline
B -- Deterministic basins & \textit{$O\!\bigl(\frac{m_{\text{ref}}}{\varepsilon}\bigr)$} & $O\!\bigl(N\,T_{\text{det}}\,m_{\text{ref}}\bigr)$ \\ \hline
C -- Random basins (two-pass) & $O\!\bigl(m_{\text{ref}}N\,T_{\text{rnd}}\bigr)$ & $O\!\bigl(m_{\text{ref}}(N\,T_{\text{rnd}} + N\,T_{\text{dyn}})\bigr)$ \\ \hline
\end{tabular}
\caption{Asymptotic time complexity of the three GPU-accelerated routines.  
For Algorithm A the best and worst cases coincide, so a single bound spans both columns.}
\end{table}

\newpage

\section{Critical Noise Analysis}
\label{sec:critical_noise} 

\subsection{Introduction and definition}
The critical noise is a noise value of interest - we notice that for no noise, the Domenicali map has precisely two well defined attractors, and from empirical results, the same map has only one attractor for noise values \( 8\times10^{-3}\) and above. Thus, there must exist a noise value, the critical noise, at which the bifurcation occurs.

\subsection{Upper Bound via Algorithm A}
Figure~\ref{fig:density-sigma0500} shows high-resolution density fields computed on an $800\times800$ grid (with $N=50\,000$ points and $T=50\,000$ steps).  At $\sigma= 4.4\times10^{-3}$ the inner filament and outer ring first make contact.  Because attractor collision cannot occur after coalescence of the two density ridges, we take this as a visual \emph{upper bound} on the critical noise $\sigma_c$.

\subsection{Upper Bound via Algorithm C}
Using the deterministic basin labels from Algorithm B [\ref{subsec:algB}], we evolve \(N = 4\times10^6\) grid points under \(T = 10\,000\) noisy iterations.  Each trajectory is colored by its original basin (red or black) and switches to white upon first escape.  

Figure~\ref{fig:random_basin_verified} shows examples at two low-noise levels:
\begin{itemize}
  \item \(\sigma = 5\times10^{-4}\): almost no escapes occur (few or no white points).
  \item \(\sigma =2\times10^{-3}\): initial white points appear along the basin boundary.
\end{itemize}

\begin{figure}[H]
  \centering
  \begin{subfigure}[b]{0.45\textwidth}
    \centering
    \includegraphics[width=0.9\linewidth]{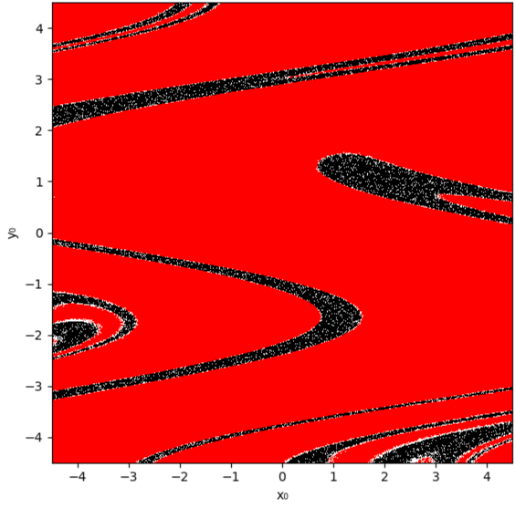}
    \caption{Basin-escape snapshot at \(\sigma=4.5\times10^{-3}\). Widespread white points show that most outer-basin trajectories escape within \(100{,}000\) steps.}
    \label{fig:basinnoise}
  \end{subfigure}
  \hfill
  \begin{subfigure}[b]{0.45\textwidth}
    \centering
    \includegraphics[width=0.9\linewidth]{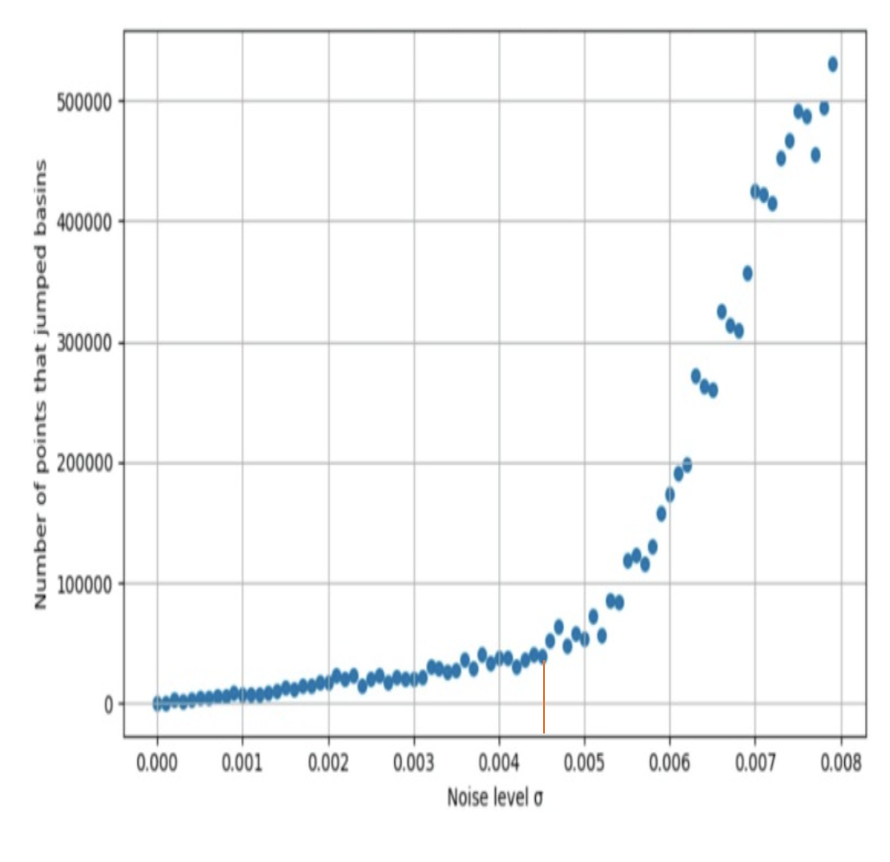}
    \caption{Spatial distribution of escaped (white) points at \(\sigma=4.5\times10^{-3}\).}
    \label{fig:plot_of_white_points}
  \end{subfigure}
  \label{fig:basinnoise_and_white}
\end{figure}

To capture the collision onset more precisely, we add Figure~\ref{fig:basinnoise}, which displays the same basin-escape plot at \(\sigma = 4.5\times10^{-3}\).  At this noise level the outer basin is densely peppered with white, indicating widespread escapes into the inner region. Given this, we take $4.5\times10^{-3}$ as an upper bound on the critical noise \(\sigma_c\). 

Figure \ref{fig:plot_of_white_points} is a graph of verified white points vs noise level at \(100,000\) iterations. We see a steady and slow increase in white points at the beginning, likely due to numerical errors (theoretically it should be exactly 0 up to the critical noise). Around $4.0\times10^{-3}$ to $5.0\times10^{-3}$, the rate of increase in white points rapidly increases to a new steeper, steady rate.

\subsection{Algorithm D: A Three-Step Algorithm for estimating the Critical Noise }
\label{sub:algD}


\subsubsection{Overview}
We now describe a numerical procedure, one introduced for the first time now, to compute an estimate of the critical noise level \(\sigma_c\).

We work on a uniform grid of resolution $N\times N$ covering the phase-space box $[-4.5,4.5]\times[4.5,4.5]$.  Our goal is to compute a noise threshold $\sigma_c$ such that for any $\sigma\ge\sigma_c$, the \emph{noisy attractor} visits at least one ``exit'' cell of the deterministic basin.  The algorithm has three main steps, but the idea is one. A basin switch must occur in \emph{one step}, because we are under a discrete mapping. So, we consider points in a fine grid, and at what noise they switch basins. Then, we investigate when a point with sufficiently low switching noise is "covered" by a random attractor at noise $\sigma$ for the first $\sigma$. This $\sigma$ forms an upper bound for $\sigma_c$

\begin{figure}[H]
  \centering
  \begin{subfigure}[b]{0.48\textwidth}
    \centering
    \includegraphics[width=\linewidth]{figures/Figure_7.png}
    \caption{No noise (\(\sigma = 0\)).}
    \label{fig:basin-no-noise}
  \end{subfigure}
  \hfill
  \begin{subfigure}[b]{0.48\textwidth}
    \centering
    \includegraphics[width=\linewidth]{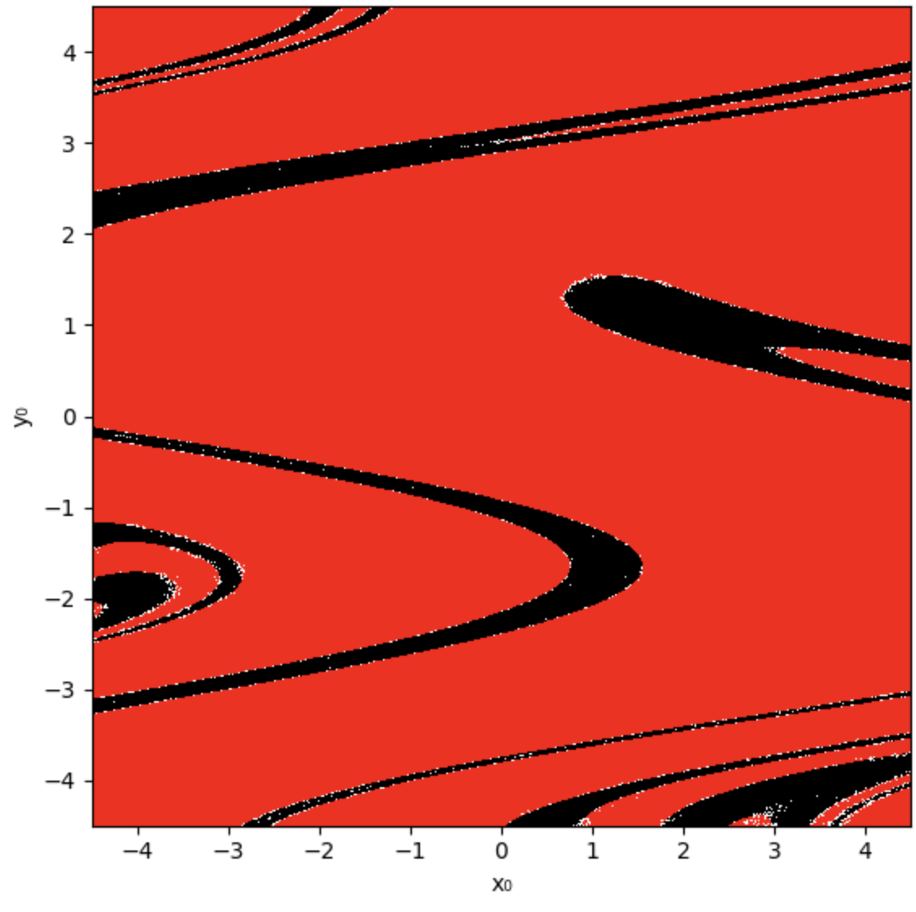}
    \caption{High noise (\(\sigma = 4\times10^{-3}\)).}
    \label{fig:basin-sigma-0.0040}
  \end{subfigure}

  \caption{The basins look similar, with the right basin (with noise) having a thicker boundary}
  \label{fig:hist_density2}
\end{figure}

We stress that Algorithm D yields an \emph{estimate} \(\hat\sigma_c\) of the true critical noise, not merely a one-sided bound. By computing the minimal noise at which any grid cell on the \emph{deterministic} basin boundary exits in one step, we obtain a consistent estimator of \(\sigma_c\): as the grid is refined and noise levels swept more finely, \(\hat\sigma_c\) converges to the true bifurcation point. Although the deterministic boundary is used for tractability, the \emph{random} basin boundary ``fattens'' under noise and quickly aligns with the deterministic one for \(\sigma\lesssim4\times10^{-3}\) (see Figure~\ref{fig:hist_density2}), so the estimate remains accurate in the regime of interest.

\subsubsection{Step 1: Compute Deterministic Basin Boundary.}
\emph{Input:} deterministic map $f(x,y;a,b)$, parameters $(a_0,b_0)$.\\
\emph{Output:} Boolean array $\mathrm{boundary}[i,j]$ indicating the basin boundary.

\begin{enumerate}
  \item[\textbf{A.}] Follow Algorithm B [\ref{subsec:algB}] to get a 2d array containing the deterministic basin.
  \item[\textbf{B.}] For each cell $(i,j)$, examine its four von Neumann neighbors $(i\pm1,j)$, $(i,j\pm1)$.  If not all labels agree, set
    \[
      \mathrm{boundary}[i,j] \;=\; \text{True},
    \]
    otherwise $\mathrm{boundary}[i,j]=\text{False}$.
\end{enumerate}
\begin{figure}[H]
  \centering
  \includegraphics[width=0.7\textwidth]{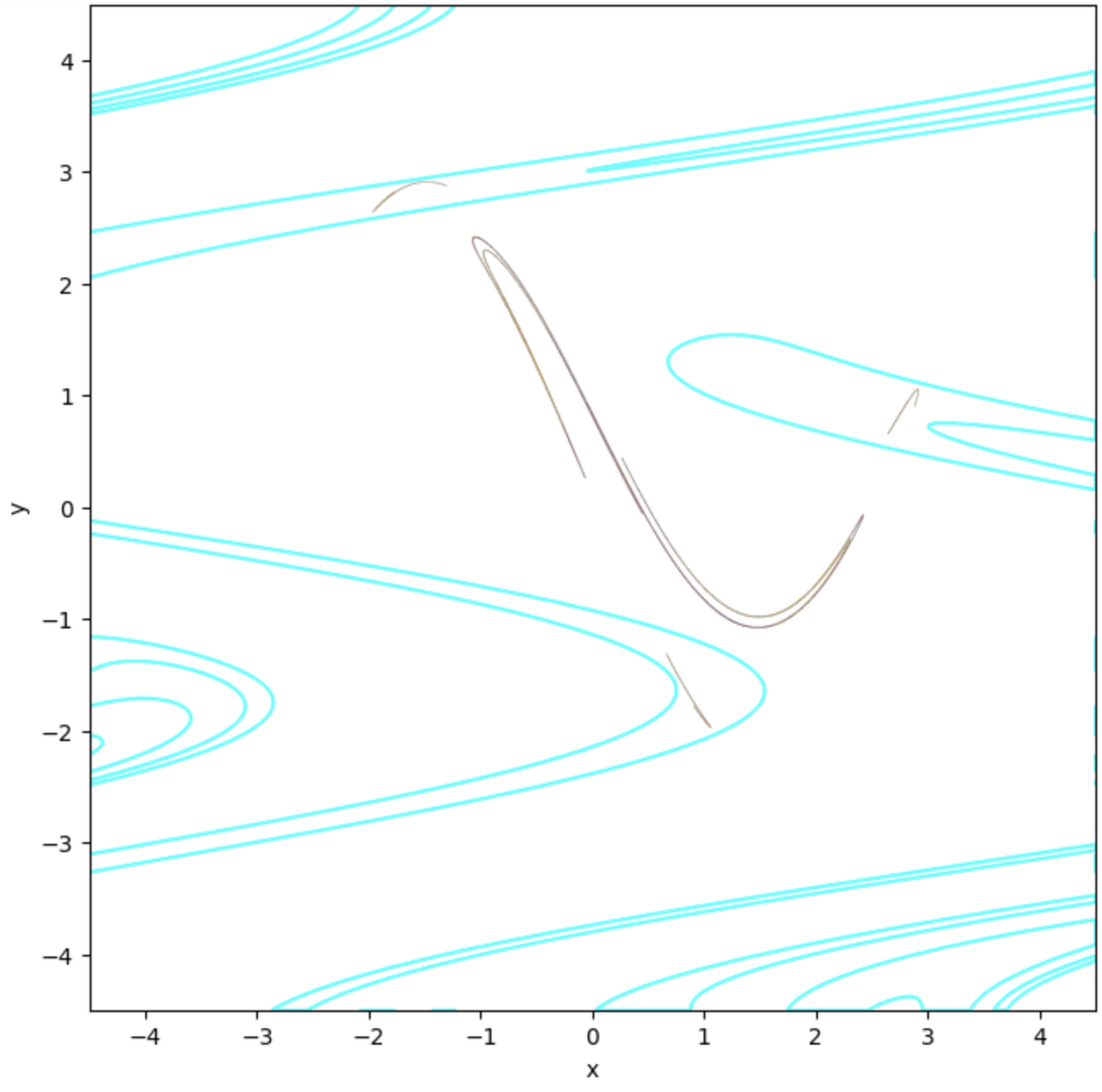}
  \caption{Basin boundary produced at step 1, overlaid with deterministic basin. Note how tightly the outer attractor fits into its basin.}
  \label{fig:basinboundary}
\end{figure}

\subsubsection{Step 2: Minimal Exit-Noise Labeling.}
\emph{Input:} boundary mask $\mathrm{boundary}[i,j]$, random map $f_\sigma$; noise range $\sigma\in[0,\sigma_{\max}]$.\\
\emph{Output:} array $\varepsilon_{\min}[i,j]$ of minimal exit noise (or $-1$ if $> \sigma_{\max}$).

\begin{enumerate}
  \item For each grid-point $(i,j)$:
    \begin{itemize}
      \item Initialise $\varepsilon_{\min}[i,j]\leftarrow +\infty$.
      \item For a finite set of noise levels $\{\sigma_k\}_{k=1}^K\subset[0,\sigma_{\max}]$ (e.g.\ linear grid):
        \begin{enumerate}
          \item Draw random parameters $a=a_0 + \delta a$, $b=b_0 + \delta b$ with $\delta a,\delta b\sim\mathcal{U}(-\sigma_k,\sigma_k)$.
          \item Compute $(x',y') = f(x_i,y_j;a,b)$.  If $(x',y')$ lies on or beyond the deterministic boundary:
            $$\mathrm{boundary}[\,\mathrm{bin}(x',y')\,]=\text{True},$$
            then set
            $$\varepsilon_{\min}[i,j] \;\leftarrow\; \min\bigl(\varepsilon_{\min}[i,j],\;\sigma_k\bigr).$$
        \end{enumerate}
      \item If $\varepsilon_{\min}[i,j]>\sigma_{\max}$, reset $\varepsilon_{\min}[i,j]\leftarrow-1$.
    \end{itemize}
\end{enumerate}

\begin{figure}[H]
  \centering
  \caption*{\textbf{Heat map of minimal exit noise}}
  \includegraphics[width=0.8\linewidth]{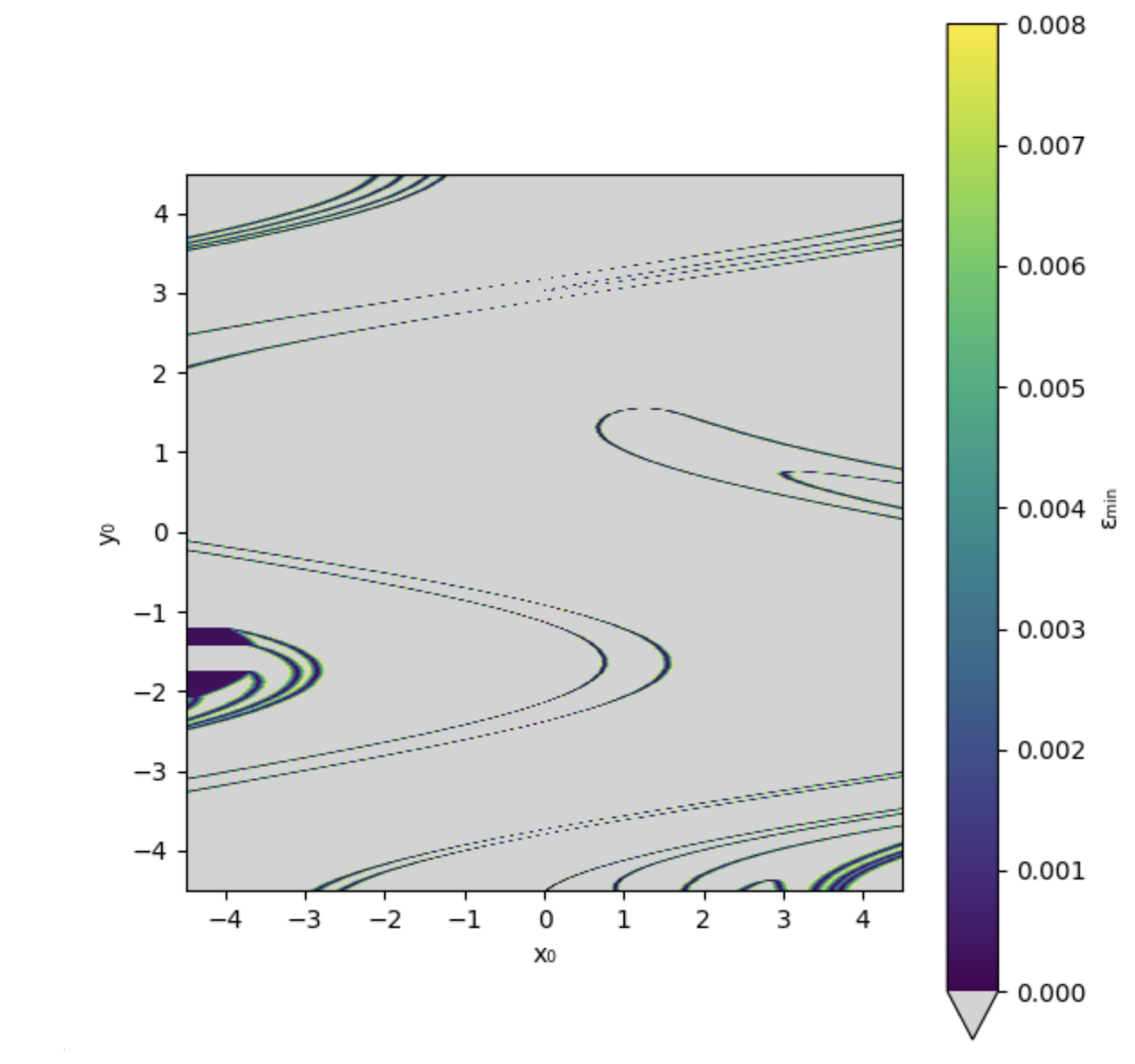}
  \caption{We set max noise to $8\times10^{-3}$. The result is then mapped to a colour map, showing that only the boundary has cells with sufficiently low exit noise. Note that  grey points require more than $8\times10^{-3}$ noise}
  \label{fig:exitnoiseheat}
\end{figure}

\subsubsection{Step 3: Noisy-Attractor Sweeping.}
\emph{Input:} minimal-noise map $\varepsilon_{\min}[i,j]$, noisy map $f_\sigma$, sweep values $\{\sigma_k\}$.\
\emph{Output:} first-hit noise $\sigma_c$.

\begin{enumerate}
  \item For each $\sigma_k$ in ascending order:
    \begin{itemize}
      \item Simulate $N$ walkers under $f_\sigma$ for $T$ steps (with burn-in), binning positions to grid cells $(i,j)$ at intervals.
      \item If any walker lands on a cell $(i,j)$ with $0\le\varepsilon_{\min}[i,j]\le\sigma_k$, then
        \[
          \sigma_c \geq \sigma_k,
        \]
        and terminate the loop.
    \end{itemize}
  \item If none found up to $\sigma_{\max}$, report ``no escape.''
\end{enumerate}

\vspace{1em}
\subsubsection{Algorithm D Conclusion}  
By construction, the first \(\sigma_k\) at which the noisy attractor hits a cell with \(\varepsilon_{\min}\le\sigma_k\) provides a natural \emph{estimate} \(\hat\sigma_c\) of the true critical noise:
\[
  \hat\sigma_c \;=\; \min\bigl\{\sigma_k : \exists\,(i,j)\text{ with } \varepsilon_{\min}[i,j]\le\sigma_k\bigr\}.
\]
In our repeated experiments, this procedure produced estimates in the range \([2.3\times10^{-3},\,2.9\times10^{-3}]\), with the finest resolution giving \(\hat\sigma_c \approx 2.3\times10^{-3}\).
Moreover, the computational complexity of Algorithm D has time complexity \textit{$O\!\bigl(\frac{m_{\text{ref}}}{\varepsilon} + Nr + NTK\bigr)$}, where N is the number of points, $m_{\text{ref}}$ is the number of points used in the reference cloud, $\varepsilon$ is the threshold when computing basins, r is the number of realisations tested per point during step 2, T is the number of iterations, and K is the number of noise values swept.

\subsection{Comparing results}

Table~\ref{tab:compare-sigma} summarizes the values of the critical noise level \(\sigma_c\)
obtained by Algorithms A, C, and D.

\begin{table}[ht]
  \centering
  \begin{tabular}{lll}
    \toprule
    Method         & Value for \(\sigma_c\)                & Notes \\
    \midrule
    Algorithm A    & \(\sigma_c \lesssim 4.4\times10^{-3}\) & Attractors contact in density field \\
    Algorithm C    & \(\sigma_c \lesssim 4.5\times10^{-3}\) & First widespread escapes in basin-escape plot \\
    Algorithm D    & \(\hat\sigma_c \approx 2.3\times10^{-3}\)      & Consistent estimate from local exit-noise sweep \\
    \bottomrule
  \end{tabular}
  \caption{Comparison of upper bounds (Algorithms A and C) and the estimate (Algorithm D) for the critical noise \(\sigma_c\).}
  \label{tab:compare-sigma}
\end{table}

Algorithms A [\ref{subsec:algA}] and C [\ref{subsec:algC}] provide very similar \emph{upper bounds} on \(\sigma_c\) (approximately \(4.4\)--\(4.5\times10^{-3}\)), based on visual or statistical detection of attractor collision and basin escape. Algorithm D [\ref{sub:algD}] , by contrast, yields a direct \emph{estimate} \(\hat\sigma_c \approx 2.3\times10^{-3}\) by making use of minimal exit noises. Since this estimate lies well within the upper-bound range, we conclude that the true critical noise satisfies
\[
  \hat\sigma_c \;\approx\; 2.3\times10^{-3}
  \quad\text{and}\quad
  \sigma_c \;\le\; 4.4\times10^{-3}.
\]





\newpage

\section{Conclusion and Future Work}
\label{sec:conclusion}

\subsection{Summary}
In this paper we explored how small, bounded fluctuations in map parameters reshape both the geometry and stability of chaotic attractors in a two-dimensional system. While we focused on the Domenicali map, using uniformly drawn coefficients at each iteration, as a concrete example, our finite-time probabilistic framework and central limit theorem for the maximal Lyapunov exponent apply to \emph{any} discrete 2D map under i.i.d.\ parameter noise. This result provides explicit formulas for the exponent's mean, variance, and finite-sample bias at arbitrary noise amplitude and over any finite observation window. 

Building on that foundation, we showed how parameter uncertainty alters the full Lyapunov spectrum and resulting fractal dimension: at very low noise levels, a subtle dip appears driven by the curvature of the deterministic spectrum, and as uncertainty increases, the gap between expansion and contraction widens systematically. On the computational side, we introduced three GPU-accelerated algorithms for separating coexisting attractors under noise, high-resolution density histograms for geometric deformation, reference-cloud classification to assign points to inner or outer attractors, and a robust two-pass verification to label basin-switching trajectories. 

Finally, by combining visual coalescence in density fields, long-term escape statistics, and a novel minimum-escape-time grid sweep, we pinpointed an estimate for the critical noise of \(\sigma\approx2.3\times10^{-3}\). These theoretical insights and practical tools now enable accurate, real-time assessment of chaos and basin stability in noisy discrete-time systems far beyond the Domenicali example.

\subsection{Future Work}

To extend our work further, four directions stand out:

\subsubsection{Lyapunov Exponent Homogenization} Investigate why the finite-time Lyapunov exponents of the inner and outer attractors remain nearly identical under noise.  Possible approaches include:
\begin{itemize}
    \item comparing the conditional distributions of $\ln\|Df\|$ sampled on each attractor
    \item analysing the repeller's unstable manifold to identify shared dominant expansion directions
    \item verifying the phenomenon at other $(a,b)$ parameter values to determine whether this matching of exponents is coincidental or a robust feature.
\end{itemize}

\subsubsection{State-Noise vs. Parameter-Noise Comparison}
    Real systems suffer both state disturbances and parameter perturbation.  A systematic study would compare
    \[
      x_{t+1} = f(x_t,y_t) + \xi_t,\quad \xi_t\sim\mathcal{N}(0,\sigma^2),
    \]
    against
    \[
      (a_t,b_t)\sim\mathcal{U}[\mu-\delta,\mu+\delta],\quad
      (x_{t+1},y_{t+1}) = f_{a_t,b_t}(x_t,y_t),
    \]
    examining shifts in the mean and variance of \(\lambda_T\), changes in Kaplan--Yorke dimension, and misclassification rates of coexisting regimes.

\subsubsection{Temporally Correlated Parameter Noise}
    In many applications parameter drift exhibits memory. Instead of analysing i.i.d.\ coefficients, one can look at how adding correlation between them changes the system.
    For example, we could define:
    \[
      a_{t+1} = \alpha\,a_t + (1-\alpha)\,\mu + \eta_t,
      \qquad
      \eta_t \sim \mathcal{N}(0,\sigma^2),
    \]
    where \(\mu\) is the long-run mean, \(\alpha\in[0,1)\) controls persistence (\(\alpha=0\) recovers the i.i.d.\ case; \(\alpha\to1\) yields strong correlation), and \(\eta_t\) is the noise.  One would study how \(\alpha\) affects FTLE-CLT convergence, the width of the chaos--period coexistence band, and the robustness of our repeller-based partitioning algorithm.

\subsubsection{Higher-Dimensional Maps under Parameter Noise}
    Many real-world problems require \(n>2\), where hyperchaos and multiple positive exponents appear.  We propose extending our finite-time theory and separation algorithm to three-dimensional (or higher) polynomial maps with i.i.d.\ parameter noise.  Key questions include:
    \begin{itemize}
      \item How does the full Lyapunov spectrum behave under noise?
      \item How do Kaplan--Yorke dimension estimates scale with dimension?
      \item Does our FTLE-histogram plus repeller-manifold refinement still reliably partition phase space as \(n\) grows?
    \end{itemize}

\newpage

\appendix

\section{Probabilistic Tools}
\label{app:prob}

\subsection{Uniform Distribution PDF and CDF}
\label{app:pdf}
Here we collect the standard formulas for 
\(\displaystyle X\sim\mathcal{U}[\mu-\delta,\mu+\delta]\).

The probability density function is
\[
  f_X(x)
  = \begin{cases}
      \dfrac{1}{2\delta}, & \mu-\delta \le x \le \mu+\delta,\\[6pt]
      0,                   & \text{otherwise},
    \end{cases}
\]
and the cumulative distribution function is
\[
  F_X(x)
  = \int_{-\infty}^{x} f_X(t)\,\mathrm{d}t
  = \begin{cases}
      0,                                         & x < \mu-\delta,\\[4pt]
      \dfrac{x - (\mu-\delta)}{2\delta},         & \mu-\delta \le x \le \mu+\delta,\\[6pt]
      1,                                         & x > \mu+\delta.
    \end{cases}
\]

\begin{figure}[ht]
  \centering
  \includegraphics[width=0.45\textwidth]{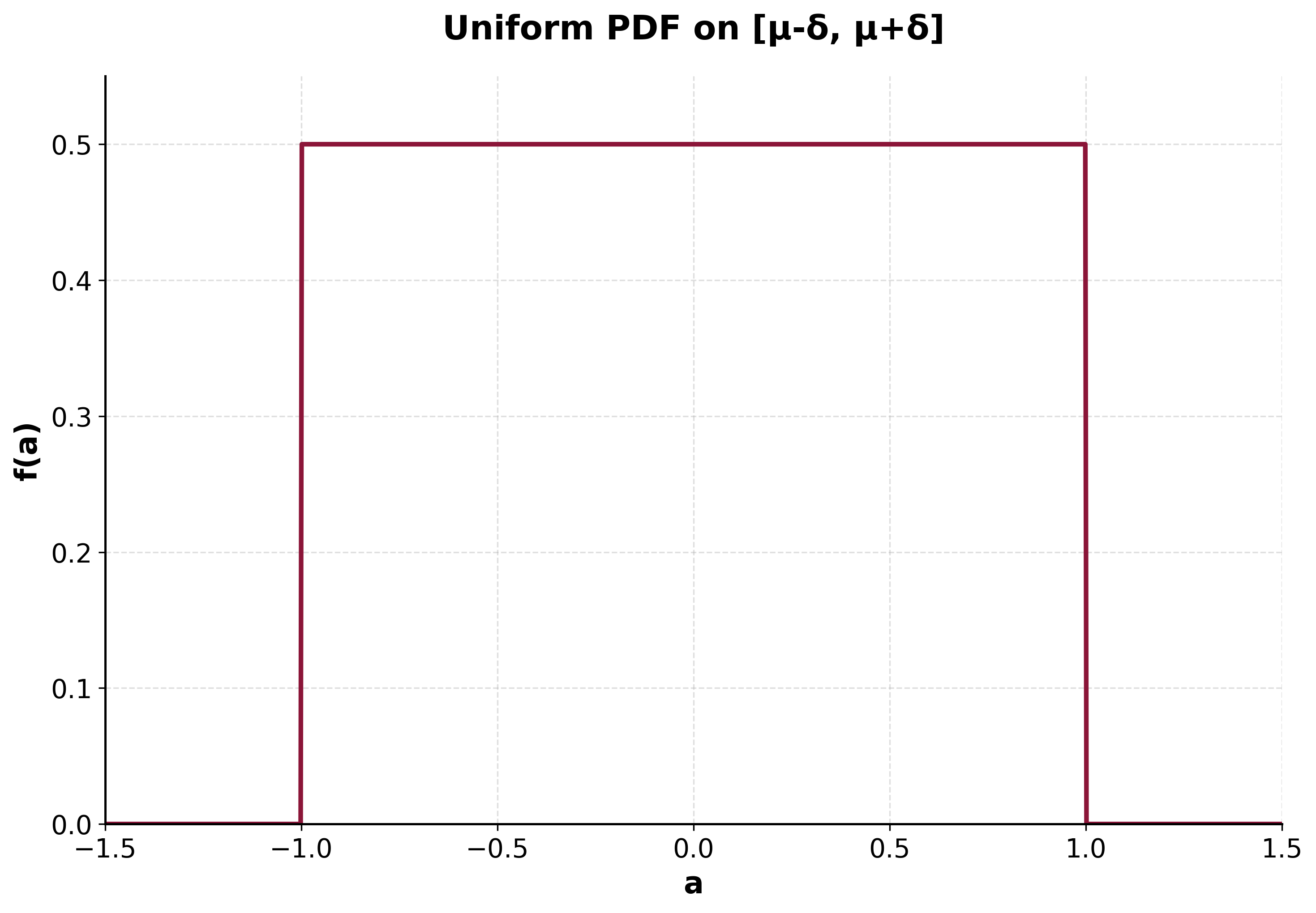}
  \hfill
  \includegraphics[width=0.45\textwidth]{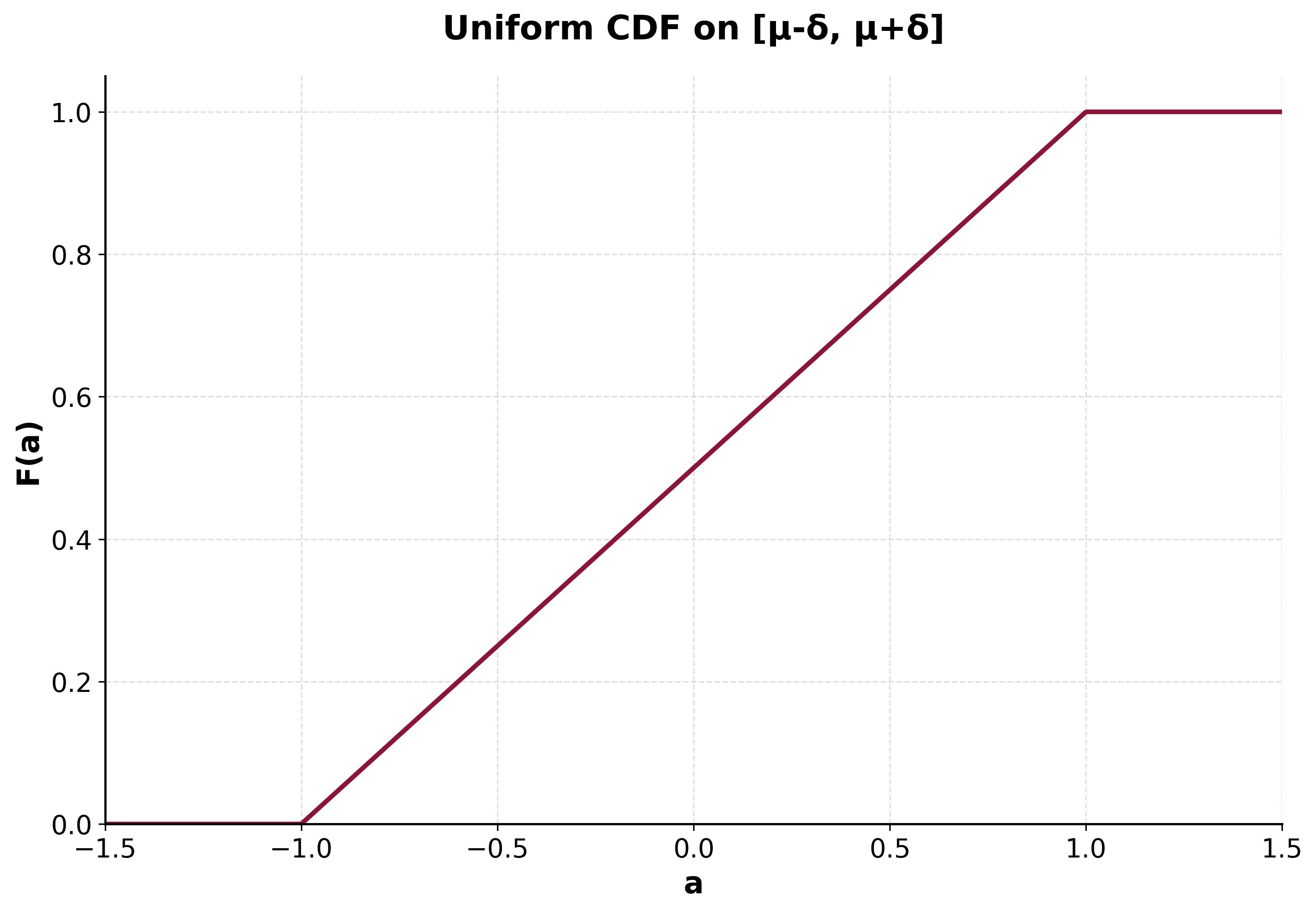}
  \caption{Left: PDF of \(\mathcal{U}[\mu-\delta,\mu+\delta]\).  Right: CDF of \(\mathcal{U}[\mu-\delta,\mu+\delta]\).}
  \label{fig:uniform-dist}
\end{figure}

\subsection{Modes of Convergence}
\label{app:convergence}

\subsubsection{Convergence in probability}
A sequence of random variables \(X_T\) converges in probability to \(X\), written
\[
  X_T \xrightarrow{p} X,
\]
if for every \(\epsilon>0\),
\[
  \lim_{T\to\infty} \Pr\bigl(|X_T - X| > \epsilon\bigr) = 0.
\]

\subsubsection{Convergence in distribution}
\(X_T\) converges in distribution to \(X\), written
\[
  X_T \xrightarrow{d} X,
\]
if for every continuity point \(x\) of the cumulative distribution function \(F_X\),
\[
  \lim_{T\to\infty} F_{X_T}(x) \;=\; F_X(x).
\]

\subsection{Oseledets--LaPage Central-Limit Theorem}
\label{app:oseledets}
Let \(\{A_t\}_{t\ge1}\) be i.i.d.\ invertible random matrices in \(\mathbb{R}^{d\times d}\) with finite exponential moments and a simple top Lyapunov exponent \(\mu\).  Then there exists \(\sigma^2>0\) such that
\[
  \sqrt{T}\,\Bigl(\tfrac1T\ln\|A_T\cdots A_1\| - \mu\Bigr)
  \;\xrightarrow{d}\;
  \mathcal{N}(0,\sigma^2).
\]
Key hypotheses include exponential integrability and strong irreducibility of the matrix law.

\subsection{Classical Central-Limit Theorem}
\label{app:classicalclt}
If \(X_i\) are i.i.d.\ with \(\mathbb{E}[X_i]=\mu\), \(\mathrm{Var}(X_i)=\sigma^2<\infty\), then
  \(\frac1{\sqrt{T}}\sum_{i=1}^T(X_i-\mu)\xrightarrow{d}\mathcal{N}(0,\sigma^2).\)

\subsection{Martingales and Martingale Difference Sequences}
\label{app:martingale}

Let \((\Omega,\mathcal F,\{\mathcal F_t\},\Pr)\) be a filtered probability space.  A sequence of integrable random variables \(\{M_t\}_{t\ge0}\) is a \emph{martingale} with respect to the filtration \(\{\mathcal F_t\}\) if for all \(t\ge0\):
\begin{enumerate}
    \item \(M_t\) is \(\mathcal F_t\)--measurable.
    \item \(\mathbb{E}[\,|M_t|\,]<\infty.\)
    \item \(\mathbb{E}[\,M_{t+1}\mid\mathcal F_t\,] = M_t.\)
\end{enumerate}

The \emph{martingale differences} are \(D_t = M_t - M_{t-1}\), which satisfy \(\mathbb{E}[D_t\mid\mathcal F_{t-1}]=0\).

\subsection{Martingale Central-Limit Theorem}
\label{app:martclt}
Let \(\{M_{T,k}\}\) be a martingale [\ref{app:martingale}] with differences \(D_{T,k}=M_{T,k}-M_{T,k-1}\) satisfying
\[
  \sum_{k=1}^T \mathbb{E}[D_{T,k}^2\mid\mathcal{F}_{k-1}]\xrightarrow{p}\sigma^2,
  \quad
  \sum_{k=1}^T \mathbb{E}[D_{T,k}^2\mathbf{1}_{\{|D_{T,k}|>\epsilon\}}\mid\mathcal{F}_{k-1}]\xrightarrow{p}0.
\]
Then \(M_{T,T}\xrightarrow{d}\mathcal{N}(0,\sigma^2)\).

\subsection{Slutsky's Lemma}
\label{app:slutsky}
If \(Y_T\xrightarrow{d}Y\) and \(Z_T\xrightarrow{p}c\), then
  \(Y_T + Z_T\xrightarrow{d}Y+c\) and \(Y_T\,Z_T\xrightarrow{d}c\,Y\).

\subsection{Delta Method and Asymptotic Normality}
\label{app:delta}

\subsubsection{Univariate Delta Method}
Suppose 
\[
  \sqrt{T}\,(\hat\theta_T - \theta)\;\xrightarrow{d}\;\mathcal{N}(0,\tau^2),
\]
and let \(g:\mathbb{R}\to\mathbb{R}\) be continuously differentiable at \(\theta\) with \(g'(\theta)\neq0\).  Then
\[
  \sqrt{T}\,\bigl(g(\hat\theta_T)-g(\theta)\bigr)
  \;\xrightarrow{d}\;
  \mathcal{N}\!\Bigl(0,\,[g'(\theta)]^2\,\tau^2\Bigr).
\]
Moreover, a second-order Taylor expansion
\[
  g(\hat\theta_T)
  = g(\theta)
    + g'(\theta)\,(\hat\theta_T-\theta)
    + \tfrac12\,g''(\theta)\,(\hat\theta_T-\theta)^2
    + o_p(T^{-1})
\]
implies a finite-sample bias
\(\mathbb{E}[g(\hat\theta_T)] - g(\theta)\approx \tfrac{g''(\theta)\,\tau^2}{2T}.\)

\subsubsection{Multivariate Delta Method}
\label{app:mdelta}

Let \(\hat\theta_T \in \mathbb{R}^d\) satisfy
\[
  \sqrt{T}\,\bigl(\hat\theta_T - \theta\bigr)
  \;\xrightarrow{d}\;
  \mathcal{N}_d\!\bigl(0,\Sigma\bigr),
\]
and let \(g:\mathbb{R}^d\to\mathbb{R}^k\) be continuously differentiable at \(\theta\).  Denote the Jacobian at \(\theta\) by
\[
  J = Dg(\theta)
    = \begin{pmatrix}
        \dfrac{\partial g_1}{\partial \theta_1} & \cdots & \dfrac{\partial g_1}{\partial \theta_d} \\[6pt]
        \vdots & \ddots & \vdots \\[6pt]
        \dfrac{\partial g_k}{\partial \theta_1} & \cdots & \dfrac{\partial g_k}{\partial \theta_d}
      \end{pmatrix}.
\]
Then
\[
  \sqrt{T}\,\bigl(g(\hat\theta_T)-g(\theta)\bigr)
  \;\xrightarrow{d}\;
  \mathcal{N}_k\!\bigl(0,\,J\,\Sigma\,J^\top\bigr).
\]
In particular, the \(i\)th component \(g_i(\hat\theta_T)\) is asymptotically normal with variance \(\bigl(J\,\Sigma\,J^\top\bigr)_{ii}\).

\subsection{Shapiro--Wilk Test}
\label{app:shapirowilk}
For samples \(\{Z_i\}\), the Shapiro--Wilk statistic
\[
  W = \frac{\bigl(\sum_i a_i Z_{(i)}\bigr)^2}{\sum_i (Z_i-\bar Z)^2},
\]
tests normality; extremely large samples can yield significant \(p\)-values even for minor deviations.

\subsection{Jensen's Inequality}
\label{app:jensen}

Let \(\varphi:\mathbb{R}\to\mathbb{R}\) be a convex function and \(X\) a real random variable with \(\mathbb{E}[|X|]<\infty\) and \(\mathbb{E}[|\varphi(X)|]<\infty\).  Then
\[
  \varphi\!\bigl(\mathbb{E}[X]\bigr)
  \;\le\;
  \mathbb{E}\!\bigl[\varphi(X)\bigr].
\]
In fact, for any probability mass (or density) function \(p\), one has
\[
  \varphi\!\Bigl(\sum_i p_i x_i\Bigr)
  \;\le\;
  \sum_i p_i\,\varphi(x_i),
  \quad
  \varphi\!\Bigl(\int x\,p(x)\,\mathrm{d}x\Bigr)
  \;\le\;
  \int \varphi(x)\,p(x)\,\mathrm{d}x.
\]

\section{Stochastic-Process}
\label{app:sp}

\subsection{Markov-Chain Basics}
\label{app:markov}
A sequence \(\{X_t\}\) with transition kernel \(P(x,A)=\Pr(X_{t+1}\in A\mid X_t=x)\) is irreducible, aperiodic, and ergodic if it admits a unique stationary distribution \(\pi\) with \(P^t(x,\cdot)\to\pi\).

\subsection{Application to Jacobian-Norm Products}
Viewing \((x_t,y_t)\to(x_{t+1},y_{t+1})\) as a Markov chain on a trapping region, the norms \(\|Df(x_t,y_t)\|\) form a function of this chain.  Under ergodicity and finite second moments, Markov-chain CLTs apply.

\section{Analytic and Linear-Algebra Tools}
\label{app:linalg}

\subsection{Spectral Operator Norm}
\label{app:specnorm}
For any \(M\in\mathbb{R}^{n\times n}\),
\[
    \|M\| =
  \|M\|_2
  = \sup_{\|v\|_2=1}\|M\,v\|_2
  = \sigma_{\max}(M)
  = \sqrt{\lambda_{\max}(M^\top M)}.
\]

\subsection{Strong Irreducibility and Proximality (SIP)}
\label{app:sip}

\subsubsection{Strong Irreducibility}
A family $\mathcal A\subset GL(d,\mathbb{R})$ is \emph{strongly irreducible} if there is no finite collection of proper subspaces $V_1,\dots,V_k\subsetneq\R^d$ whose union is invariant under every $A\in\mathcal A$.

\subsubsection{Proximality}
A matrix $A\in GL(d,\mathbb{R})$ is \emph{proximal} if it has a unique real eigenvalue of maximal modulus, with one-dimensional eigenspace.  A family $\mathcal A$ is \emph{proximal} if its semigroup contains at least one proximal element.

\subsection{Taylor Expansion and Hessian Bounds}
\label{app:taylor}
For a smooth \(h:\mathbb{R}^2\to\mathbb{R}\) near \((x^*,y^*)\),
\[
  h(x,y)
  = h(x^*,y^*)
    + \nabla h(x^*,y^*)^\top\!\begin{pmatrix}x-x^*\\y-y^*\end{pmatrix}
    + \tfrac12\,(x-x^*,y-y^*)\,H\,(x-x^*,y-y^*)^\top
    + O(\|(x,y)-(x^*,y^*)\|^3),
\]
with Hessian \(H\).

\subsection{Bias in Finite-Time Lyapunov Estimates}
\label{app:bias}
Expanding \(\ln\|Df\|\) via [\ref{app:specnorm}] and taking expectations gives a bias 
\(\mathbb{E}[\lambda_T]-\lambda_\infty\approx \tfrac1{2T}\mathrm{trace}(H\,\Sigma)\), 
where \(\Sigma\) is the invariant-measure covariance.

\subsection{Proof of the Difference-of-Normals Statement}
\label{app:diffnorm}
If \(Z_1\sim\mathcal{N}(\mu_1,\sigma_1^2)\) and \(Z_2\sim\mathcal{N}(\mu_2,\sigma_2^2)\) independent, then
\[
  Z_1 - Z_2 \sim \mathcal{N}(\mu_1 - \mu_2,\;\sigma_1^2 + \sigma_2^2).
\]

\begin{proof}
Let \(Z_1\sim\mathcal{N}(\mu_1,\sigma_1^2)\) and \(Z_2\sim\mathcal{N}(\mu_2,\sigma_2^2)\) be independent.  We compute the characteristic function of \(W = Z_1 - Z_2\):
\[
\phi_W(t)
= \mathbb{E}\bigl[e^{i t (Z_1 - Z_2)}\bigr]
= \mathbb{E}\bigl[e^{i t Z_1}\bigr]\;\mathbb{E}\bigl[e^{-i t Z_2}\bigr]
= \phi_{Z_1}(t)\,\phi_{Z_2}(-t).
\]
Since for a normal \(Z\sim\mathcal{N}(\mu,\sigma^2)\) one has
\[
\phi_Z(t)
= \exp\!\Bigl(i\mu t \;-\;\tfrac12\sigma^2 t^2\Bigr),
\]
we get
\[
\phi_W(t)
= \exp\!\Bigl(i\mu_1 t - \tfrac12\sigma_1^2 t^2\Bigr)
  \;\times\;
  \exp\!\Bigl(-\,i\mu_2 t - \tfrac12\sigma_2^2 t^2\Bigr)
= \exp\!\Bigl(i(\mu_1-\mu_2)t - \tfrac12(\sigma_1^2+\sigma_2^2)t^2\Bigr).
\]
But this is exactly the characteristic function of
\(\mathcal{N}(\mu_1-\mu_2,\;\sigma_1^2+\sigma_2^2)\).  Therefore
\[
Z_1 - Z_2 \;\sim\; \mathcal{N}\bigl(\mu_1-\mu_2,\;\sigma_1^2+\sigma_2^2\bigr),
\]
as claimed.
\end{proof}

\newpage



\bibliography{bibliography}

\end{document}